# GAUSSIAN LIMITS FOR RANDOM MEASURES IN GEOMETRIC PROBABILITY

By Yu. Baryshnikov and J. E. Yukich[1]

*Bell Laboratories and Lehigh University*

We establish Gaussian limits for general measures induced by binomial and Poisson point processes in $d$-dimensional space. The limiting Gaussian field has a covariance functional which depends on the density of the point process. The general results are used to deduce central limit theorems for measures induced by random graphs (nearest neighbor, Voronoi and sphere of influence graph), random sequential packing models (ballistic deposition and spatial birth–growth models) and statistics of germ–grain models.

**1. Introduction.** The purpose of this paper is to provide a methodology for showing that renormalized random point measures in geometric probability converge weakly to a generalized Gaussian field. We focus on random point measures, defined on the Borel subsets of $\mathbb{R}^d$, of the following types:

(i)  point measures associated with random graphs in computational geometry, including nearest neighbor graphs, Voronoi graphs and sphere of influence graphs,

(ii)  point measures arising in random sequential packing models, including random sequential adsorption (RSA) and spatial birth–growth models, and

(iii)  point measures associated with germ–grain models.

The total mass of random point measures yields random functionals, which in the context of the measures (i)–(iii), have been extensively studied; see [2, 21, 23, 24, 26, 32], [5, 7, 8, 9, 11, 15, 25, 27] and [12, 13, 22, 26], respectively, as well as the references therein. With the exception of [13], the study of the random measures (i)–(iii) has received considerably less attention. We show here after renormalization that measures of the type (i)–(iii)

Received November 2002; revised January 2004.

[1]Supported in part by NSA Grant MDA904-01-1-0029 and NSF Grant DMS-02-03720.

converge to a generalized Gaussian field; that is, their finite-dimensional distributions, as described by the action of the measure on continuous test functions, converge to those of a generalized finitely additive Gaussian field. The results relate the large-scale Gaussian limit properties of renormalized random point measures to the small-scale properties of the underlying binomial or Poisson point process.

The general approach taken here, which employs stabilization of functionals and coupling arguments, has the particular benefit of describing the limiting variance over large sample sizes as a function of the underlying density of points. A similar approach is used in [26], which treats the easier problem of finding limiting means.

Random measures considered here assume the form $\sum_{x \in \mathcal{X}} \xi(x; \mathcal{X}) \delta_x$, where $\mathcal{X}$ is a random point set in $\mathbb{R}^d$, $\delta_x$ is the Dirac point measure at $x$ and $\xi(x; \mathcal{X})$ is a weight representing the interaction of $x$ with respect to $\mathcal{X}$ and is usually defined in terms of the underlying geometry. For all constants $\lambda > 0$ and probability densities $\kappa$, let $\mathcal{P}_{\lambda\kappa}$ be a Poisson point process with intensity measure $\lambda\kappa \colon \mathbb{R}^d \to \mathbb{R}^+$. Define the "binomial" point process $\mathcal{X}_n := \{X_1, \ldots, X_n\}$, where $X_i, i \geq 1$, are i.i.d. with density $\kappa$. All of our results follow from general central limit theorems (Theorems 2.1, 2.2, 2.5) which show that renormalized measures of the type

$$(1.1) \qquad \lambda^{-1/2} \sum_{x \in \mathcal{P}_{\lambda\kappa}} \xi(\lambda^{1/d} x; \lambda^{1/d} \mathcal{P}_{\lambda\kappa}) \delta_x, \qquad \lambda \geq 1,$$

as well as their respective renormalized binomial counterparts,

$$(1.2) \qquad n^{-1/2} \sum_{X_i \in \mathcal{X}_n} \xi(n^{1/d} X_i; n^{1/d} \mathcal{X}_n) \delta_{X_i}, \qquad n \geq 1,$$

converge weakly as $\lambda \to \infty$ (resp. as $n \to \infty$) to a Gaussian field with a covariance functional described in terms of the weight $\xi$ and the underlying density $\kappa$ of points.

The general central limit theorem (CLT) for the measures (1.2) implies a CLT for the "total mass" functional $\sum_{X_i \in \mathcal{X}_n} \xi(n^{1/d} X_i; n^{1/d} \mathcal{X}_n)$. $\kappa$ need not be uniform and $\xi$ need not be translation invariant, showing that even in the functional setting, we extend and generalize previous results [2, 3, 5, 24, 25].

The proofs are based on the method of cumulants, which requires showing that the cumulants of the integrals of the rescaled measures (1.1) against a large class of test functions converge to the cumulants of a normal random variable. An important tool is "stabilization" of functionals, used heavily in [5, 24, 25, 26]. Stabilization guarantees that the pair correlation function for the weights $\xi(x, \mathcal{P}_\lambda)$, $x \in \mathbb{R}^d$, decays fast enough to prove convergence of the cumulant measures associated with (1.1). To show convergence of the first- and second-order cumulant measures against test functions, we rely upon the "objective method," which exploits the fact that if $\xi$ is locally determined in a sense to be made precise, then the large $\lambda$ behavior of $\xi(\lambda^{1/d}x, \lambda^{1/d}\mathcal{P}_{\lambda\kappa})$, $x$ fixed, is approximated by the behavior of $\xi$ on homogeneous Poisson processes. This idea was developed in [26], a law of large numbers (LLN) precursor to the present paper. To show convergence of the higher-order cumulant measures, we employ cumulant expansion techniques [20].

## 2. Main results.

2.1. *Terminology.* Before stating our main results we introduce some terminology similar to that developed in [5, 24, 25, 26]. Let $\mathcal{X} \subset \mathbb{R}^d$ be finite and $y + \mathcal{X} := \{y + x : x \in \mathcal{X}\}$ for all $y \in \mathbb{R}^d$. Given $a > 0$, let $a\mathcal{X} := \{ax : x \in \mathcal{X}\}$. For $x \in \mathbb{R}^d$, $|x|$ denotes the Euclidean norm and $B_r(x)$ denotes the Euclidean ball centered at $x$ of radius $r$. $\omega_d$ denotes the volume of the unit ball in $\mathbb{R}^d$ and $\mathbf{0}$ denotes the origin of $\mathbb{R}^d$.

Throughout, $\mathcal{A}$ denotes the family of compact, convex subsets $A \subset \mathbb{R}^d$ with nonempty interior. Let $\mathcal{A}'$ denote $\mathcal{A}$ together with the space $\mathbb{R}^d$. For $A \in \mathcal{A}'$, $C(A)$ denotes the continuous functions $f : A \to \mathbb{R}$. For $f \in C(A)$ and $\mu$ a Borel measure on $\mathcal{B}(A)$ we let $\langle f, \mu \rangle := \int_A f \, du$. Given $f : \mathbb{R}^d \to \mathbb{R}$, let Supp $f$ be the closure of $\{x \in \mathbb{R}^d : f(x) \neq 0\}$.

Let $\xi(x; \mathcal{X})$ be a measurable $\mathbb{R}$-valued function defined for all pairs $(x, \mathcal{X})$, where $x$ is an element of $\mathcal{X}$. For the moment, *we assume that $\xi$ is translation invariant*, that is, $\xi(x; \mathcal{X}) = \xi(x - y; \mathcal{X} - y)$ for all $y \in \mathbb{R}^d$. When $x \notin \mathcal{X}$, we abbreviate notation and write $\xi(x; \mathcal{X})$ instead of $\xi(x; \mathcal{X} \cup x)$.

Any finite $\mathcal{X}$ induces the point measure $\sum_{x \in \mathcal{X}} \xi(x; \mathcal{X})\delta_x$. For all $\lambda > 0$, let $\xi_\lambda(x; \mathcal{X}) := \xi(\lambda^{1/d}x; \lambda^{1/d}\mathcal{X})$. A density $\kappa$ with support on $A \in \mathcal{A}$ and a weight $\xi$ generate scaled *random point measures*

$$\mu_{\lambda\kappa}^\xi := \sum_{x \in \mathcal{P}_{\lambda\kappa}} \xi_\lambda(x; \mathcal{P}_{\lambda\kappa})\delta_x.$$

The centered version of $\mu_{\lambda\kappa}^\xi$ is $\bar{\mu}_{\lambda\kappa}^\xi := \mu_{\lambda\kappa}^\xi - \mathbb{E}\mu_{\lambda\kappa}^\xi$, where for all Borel sets $B \subset A$, $\mathbb{E}[\mu_{\lambda\kappa}^\xi(B)] = \lambda \int_B \mathbb{E}[\xi_\lambda(x; \mathcal{P}_{\lambda\kappa})]\kappa(x)\,dx$. This paper develops a methodology for establishing convergence of the finite-dimensional distributions of the renormalized random point measures $\lambda^{-1/2}\bar{\mu}_{\lambda\kappa}^\xi, \lambda \geq 1$. Previous work [25, 26] developed laws of large numbers for the *total mass functional*

$$(2.1) \qquad \mu_{\lambda\kappa}^\xi(A) := \sum_{x \in \mathcal{P}_{\lambda\kappa}} \xi_\lambda(x; \mathcal{P}_{\lambda\kappa})$$

as well as CLTs [24] for translation-invariant functionals on uniform point sets which are "locally determined." The following concept of stabilization makes precise the idea of "locally determined." For all $0 \leq a < b < \infty$, let $\mathcal{F}(a, b)$ consist of all $f : \mathbb{R}^d \to \mathbb{R}^+$ having support in $\mathcal{A}'$ and such that the range of $f$ is in $[a, b] \cup \{0\}$. The common probability space $(\Omega, \mathcal{F}, P)$ for all $\mathcal{P}_f, f \in \mathcal{F}(a, b)$, can be chosen as the probability space of the Poisson point process $\mathcal{P}^*$ having intensity 1 on $\mathbb{R}^d \times \mathbb{R}^+$ such that $\mathcal{P}_f = \pi_{\mathbb{R}^d}(\mathcal{P}^* \cap \{(x, h) \in \mathbb{R}^d : h \le f(x)\})$, where $\pi_{\mathbb{R}^d}$ denotes projection from $\mathbb{R}^d \times \mathbb{R}^+$ onto $\mathbb{R}^d$. For all $\tau > 0$, let $\mathcal{P}_\tau$ denote a homogeneous Poisson point process on $\mathbb{R}^d$ with intensity $\tau$.

DEFINITION 2.1. The functional $\xi$ is *stabilizing* if for all $A \in \mathcal{A}', 0 \le a < b < \infty, \lambda > 0$, and $x \in \lambda A$, there exists an a.s. finite random variable $R(x) := R(x, \lambda, a, b, A)$ (a *radius of stabilization* for $\xi$ at $x$) defined on $(\Omega, \mathcal{F}, P)$ such that for all $f \in \mathcal{F}(a, b)$, with Supp $f = \lambda A$, and all finite $\mathcal{X} \subset \lambda A \setminus B_R(x)$ we have

$$\xi(x; (\mathcal{P}_f \cap B_R(x)) \cup \mathcal{X}) = \xi(x; \mathcal{P}_f \cap B_R(x))$$

and moreover $\sup_{x \in \mathbb{R}^d} P[R(x, \lambda, a, b, A) > t] \to 0$ as $t \to \infty$. When $\xi$ stabilizes, then for all $x \in \mathbb{R}^d$ and all $\tau > 0$ we define

$$\xi(x; \mathcal{P}_\tau) := \lim_{l \to \infty} \xi(x; \mathcal{P}_\tau \cap B_l(x)).$$

Thus $R := R(x, \lambda, a, b, A)$ is a radius of stabilization if the value of $\xi(x; \mathcal{P}_f)$, $f \in \mathcal{F}(a, b)$, is unaffected by changes outside $B_R(x)$. One might expect that exponential decay of the tails of $R$ implies exponential decay of the correlations of $\xi$ and thus convergence of $\lambda^{-1/2} \bar{\mu}_{\lambda\kappa}^\xi, \lambda \ge 1$, to a Gaussian field. This loosely formulated idea figures prominently in interacting particle systems on the lattice, and also in cluster expansions and the moment method in statistical physics [20]. Assuming neither translation invariance of $\xi$ nor spatial homogeneity of points, we will show that this idea also works well in the continuum, where it yields convergence of $\langle f, \lambda^{-1/2} \bar{\mu}_{\lambda\kappa}^\xi \rangle_\lambda, f \in C(A)$, to a Gaussian field whose covariance depends on the density of points. This motivates defining uniform tail probabilities for the radii $R(x, \lambda, a, b, A)$:

$$r(t) := r(t, a, b, A) := \sup_{x \in \mathbb{R}^d, \lambda > 0} P[R(x, \lambda, a, b, A) \ge t].$$

$r(t)$ quantifies the region of influence of points in the Poisson point sets $\mathcal{P}_f$, whenever $f \in \mathcal{F}(a, b)$ and Supp $f$ is a scalar multiple of $A$. $\xi$ is *exponentially stabilizing* if $r(t)$ decays exponentially in $t$ for all $a, b$ and any $A \in \mathcal{A}'$. $\xi$ is *polynomially stabilizing* if for all $a, b$ and $A \in \mathcal{A}'$ we have $\int_0^\infty (r(t))^{1/2} t^{d-1} \, dt < \infty$, which readily implies the rough estimate $r(t) = o(t^{-2})$.

The next condition is used frequently in the scaling limit analysis of random fields on lattices (e.g., page 193 in [20]) and it is only natural to use it in the continuum setting as well. Here and henceforth $\kappa$ is a probability density which is continuous on its support and Supp $\kappa \in \mathcal{A}$. Let $\mathcal{C}$ denote the collection of finite point sets in $\mathbb{R}^d$.

DEFINITION 2.2. $\xi$ has a moment of order $p > 0$ with respect to $\kappa$ if

$$(2.2) \qquad \sup_{\lambda > 0, x \in [0, \lambda^{1/d}]^d A, \mathcal{X} \in \mathcal{C}} \mathbb{E}[|\xi_\lambda(x; \mathcal{P}_{\lambda\kappa} \cup \mathcal{X})|^p] < \infty$$

$$\text{and for all } \lambda > 0 \qquad \sup_{x \in \mathbb{R}^d, \mathcal{X} \in \mathcal{C}} \mathbb{E}[|\xi(x; \mathcal{P}_\lambda \cup \mathcal{X})|^p] < \infty.$$

We implicitly assume for all $l > 0$ that $\xi^\ell := \xi(x, \mathcal{X} \cap B_l(x))$ has moments no larger than those of $\xi$.

2.2. *General central limit theorems.* Under stabilization and moment conditions, we will show in Theorem 2.1 that $\mathbb{E}\bar{\mu}^\xi_{\lambda\kappa}$ and $\mathrm{Var}\,\bar{\mu}^\xi_{\lambda\kappa}$ have volume order asymptotics and that the scaling limit of the finite-dimensional distributions of the renormalized measures $\lambda^{-1/2}\bar{\mu}^\xi_{\lambda\kappa}$ is a mean zero Gaussian field. Theorem 2.1 is a special case of the upcoming Theorem 2.4 and applications of both are described in Section 3.

By the convergence of finite-dimensional distributions of random signed measures $\mu_n$ to those of a generalized Gaussian field we mean the convergence in distribution of the integrals $\int f \, d\mu_n$ to the corresponding normal random variables for all test functions $f \in C(A)$. This is the usual functional analytic point of view where a measure is viewed as a continuous linear functional acting on continuous functions. Henceforth we say that *measures converge to a Gaussian field if their finite-dimensional distributions converge.*

For all $\tau > 0$, let

$$V^\xi(\tau) := \mathbb{E}[\xi^2(\mathbf{0}; \mathcal{P}_\tau)]$$
$$+ \int_{\mathbb{R}^d}(\mathbb{E}[\xi(\mathbf{0}; \mathcal{P}_\tau \cup y) \cdot \xi(y; \mathcal{P}_\tau \cup \mathbf{0})] - \mathbb{E}[\xi(\mathbf{0}; \mathcal{P}_\tau)]\mathbb{E}[\xi(y; \mathcal{P}'_\tau)])\tau \, dy,$$

where $\mathcal{P}'_\tau$ denotes an independent copy of $\mathcal{P}_\tau$.

THEOREM 2.1. (i) *If $\xi$ is stabilizing and satisfies* (2.2) *for some $p > 1$, then for all $f \in C(A)$*

$$(2.3) \qquad \lim_{\lambda \to \infty} \frac{\mathbb{E}[\langle f, \mu^\xi_{\lambda\kappa}\rangle]}{\lambda} = \int_A f(x)\mathbb{E}[\xi(\mathbf{0}; \mathcal{P}_{\kappa(x)})]\kappa(x) \, dx,$$

*whereas if $\xi$ is polynomially stabilizing and satisfies* (2.2) *for $p = 4$, then*

$$(2.4) \qquad \lim_{\lambda \to \infty} \frac{\mathrm{Var}[\langle f, \mu^\xi_{\lambda\kappa}\rangle]}{\lambda} = \int_A f^2(x)V^\xi(\kappa(x))\kappa(x) \, dx.$$

(ii) *If $\xi$ is exponentially stabilizing and satisfies* (2.2) *for all $p > 0$, then $\lambda^{-1/2}\bar{\mu}^\xi_{\lambda\kappa}$ converges as $\lambda \to \infty$ to a Gaussian field with covariance kernel $\int_A f_1(x)f_2(x)V^\xi(\kappa(x))\kappa(x) \, dx$.*

Statistical applications often require the analog of Theorem 2.1 for measures induced by exactly $n$ i.i.d. points on $A$. This "de-Poissonized" version of Theorem 2.1 goes as follows. Let $X_i$, $i \geq 1$, be i.i.d. with common density $\kappa$, $\mathcal{X}_n := \{X_1, \ldots, X_n\}$, and $\rho_n^\xi := \sum_{i=1}^n \xi_n(X_i; \mathcal{X}_n)\delta_{X_i}$ the random "de-Poissonized" measures induced by $\kappa$ and $\xi$. To obtain the convergence of the finite-dimensional distributions of $\bar{\rho}_n^\xi := \rho_n^\xi - \mathbb{E}\rho_n^\xi$ we need some additional terminology and assumptions.

Let $\mathcal{X}_{m,n}$ be a point process consisting of $m$ i.i.d. random variables $n^{1/d}X$ on $n^{1/d}A$, where $X$ has density $\kappa$. For all $\mathcal{X}$, let $H(\mathcal{X}) := H^\xi(\mathcal{X}) := \sum_{x \in \mathcal{X}} \xi(x; \mathcal{X})$ and for all $\lambda > 0$ let $H_\lambda^\xi(\mathcal{X}) := \sum_{x \in \mathcal{X}} \xi_\lambda(x; \mathcal{X})$. For any finite $\mathcal{X}$, let $\Delta_x(\mathcal{X}) := H(\mathcal{X} \cup x) - H(\mathcal{X})$. Say that $H$ satisfies the *bounded moments condition for $\kappa$* (cf. [24]) if

$$(2.5) \qquad \sup_n \sup_{x \in n^{1/d}A} \sup_{m \in [n/2, 3n/2]} \mathbb{E}[\Delta_x^4(\mathcal{X}_{m,n})] < \infty.$$

If $H$ satisfies the bounded moments condition for $\kappa$ then we will assume throughout that $H_f^\xi$ defined by $H_f^\xi(\mathcal{X}) := \sum_{x \in \mathcal{X}} f(x)\xi(x; \mathcal{X})$, $f \in C(A)$, also satisfies the bounded moments condition. This assumption is satisfied in all of our applications in Section 3.

The next definition recalls a notion of stabilization for $H$ introduced in [24]. We are grateful to Mathew Penrose for pointing out that stabilization of $H$ rather than that of $\xi$ is essential for the upcoming de-Poissonization methods of Section 6; this observation corrects an earlier version of our results.

DEFINITION 2.3.   The functional $H := H^\xi$ is *strongly stabilizing* if for all $\tau > 0$, there exist a.s. finite random variables $S$ (a *radius of stabilization* of $H$) and $\Delta^\xi(\tau)$ such that with probability 1,

$$(2.6) \qquad \Delta^\xi(\tau) = \Delta_{\mathbf{0}}((\mathcal{P}_\tau \cap B_S(\mathbf{0})) \cup \mathcal{A})$$

for all finite $\mathcal{A} \subset \mathbb{R}^d \setminus B_S(\mathbf{0})$.

If $H^\xi$ is strongly stabilizing, then we will assume throughout that $H_f^\xi$, $f \in C(A)$, is also strongly stabilizing. This assumption is satisfied in all of our applications in Section 3.

Let $D^\xi(\tau) := \mathbb{E}[\Delta^\xi(\tau)]$ for all $\tau > 0$. The following de-Poissonized version of Theorem 2.1 shows that $\mathbb{E}[\langle f, \rho_n^\xi \rangle]$ and $\mathrm{Var}[\langle f, \rho_n^\xi \rangle]$, $f \in C(A)$, have volume order fluctuations and that the scaling limit of the re-normalized measures $n^{-1/2}\bar{\rho}_n^\xi$ is a mean zero Gaussian field.

THEOREM 2.2.   (i) *If $\xi$ is stabilizing and satisfies* (2.2) *for some $p > 1$, then for all $f \in C(A)$,*

$$(2.7) \qquad \lim_{n \to \infty} \frac{\mathbb{E}[\langle f, \rho_n^\xi \rangle]}{n} = \int_A f(x)\mathbb{E}[\xi(\mathbf{0}; \mathcal{P}_{\kappa(x)})]\kappa(x)\,dx,$$

*whereas if $\xi$ is polynomially stabilizing and satisfies* (2.2) *for $p = 4$, and if $H$ is strongly stabilizing and satisfies the bounded moments condition for $\kappa$, then for all $f \in C(A)$,*

$$
\begin{aligned}
(2.8) \quad &\lim_{n \to \infty} \frac{\mathrm{Var}[\langle f, \rho_n^\xi \rangle]}{n} \\
&= \int_A f^2(x) V^\xi(\kappa(x)) \kappa(x)\, dx - \left( \int_A f(x) D^\xi(\kappa(x)) \kappa(x)\, dx \right)^2.
\end{aligned}
$$

(ii) *If $\xi$ is exponentially stabilizing and satisfies* (2.2) *for all $p > 0$ and if $H$ is strongly stabilizing, then $n^{-1/2} \bar{\rho}_n^\xi$ converges as $n \to \infty$ to a Gaussian field with covariance kernel*

$$
\begin{aligned}
(2.9) \quad &\int_A f_1(x) f_2(x) V^\xi(\kappa(x)) \kappa(x)\, dx \\
&- \int_A f_1(x) D^\xi(\kappa(x)) \kappa(x)\, dx \int_A f_2(x) D^\xi(\kappa(x)) \kappa(x)\, dx.
\end{aligned}
$$

(iii) *If the distribution of $\Delta^\xi(\kappa(X))$ is nondegenerate, then*

$$
(2.10) \qquad \lim_{n \to \infty} \frac{\mathrm{Var}[\bar{\rho}_n^\xi(A)]}{n} > 0,
$$

*that is, the limiting Gaussian field is nondegenerate.*

REMARKS. (i) Theorems 2.1 and 2.2 generalize existing central limit theorems in geometric probability (Heinrich and Molchanov [13], Malyshev [19], Penrose and Yukich [24] and Ivanoff [16]) in several ways: (a) they show asymptotic convergence of measures to a Gaussian field, thus also yielding asymptotic convergence of functionals to a limiting normal random variable, (b) they identify the limiting variance and covariance structure in terms of the underlying density of points, and (c) they do not assume spatial homogeneity of the underlying points. Theorem 2.2 implies that for $f_1, \ldots, f_m \in C(A)$, the random vector $\langle \langle f_1, n^{-1/2} \bar{\rho}_n^\xi \rangle, \ldots, \langle f_m, n^{-1/2} \bar{\rho}_n^\xi \rangle \rangle$ converges to a multivariate Gaussian random variable.

(ii) Evaluating (2.4) and (2.8) is in general difficult. However, for some problems of geometric probability, for example, those involving functionals which count the number of pairs of points within a specified distance of one another, it is relatively simple to evaluate $V^\xi$ and $D^\xi$ [6]. Moreover, a simplification of (2.4) and (2.8) occurs whenever $\xi$ is *homogeneous of order $\gamma$*, that is, whenever there is a constant $\gamma > 0$ such that $\xi$ satisfies the relation $\xi(ax; a\mathcal{X}) = a^\gamma \xi(x; \mathcal{X})$ for all positive scalars $a$ and all finite point sets $\mathcal{X} \ni x$. Homogeneity occurs naturally in many problems of geometric probability.

If $\xi$ is homogeneous of order $\gamma$, then $V^{\xi}(\tau) = V^{\xi}(1)\tau^{-2\gamma/d}$, and $D^{\xi}(\tau) = D^{\xi}(1)\tau^{-\gamma/d}$, yielding

$$(2.11) \qquad \lim_{\lambda \to \infty} \frac{\mathrm{Var}[\langle f, \bar{\mu}^{\xi}_{\lambda\kappa}\rangle]}{\lambda} = V^{\xi}(1) \int_A f^2(x)\kappa(x)^{(d-2\gamma)/d}\, dx$$

and

$$(2.12) \qquad \lim_{n \to \infty} \frac{\mathrm{Var}[\langle f, \bar{\rho}^{\xi}_{n}\rangle]}{n} = V^{\xi}(1) \int_A f^2(x)\kappa(x)^{(d-2\gamma)/d}\, dx$$
$$- (D^{\xi}(1))^2 \left(\int_A f(x)\kappa(x)^{(d-\gamma)/d}\, dx\right)^2.$$

If $\kappa$ is the uniform distribution on the unit cube, then by (2.12)

$$\lim_{n \to \infty} \frac{\mathrm{Var}[\bar{\rho}^{\xi}_{n}([0,1]^d)]}{n} = V^{\xi}(1) - (D^{\xi}(1))^2,$$

which is strictly positive whenever $\Delta^{\xi}(1)$ is nondegenerate. If $\xi$ is scale invariant, or homogeneous of order 0, then for any $\kappa$ with support $A$,

$$\lim_{n \to \infty} \frac{\mathrm{Var}[\bar{\rho}^{\xi}_{n}(A)]}{n} = V^{\xi}(1) - (D^{\xi}(1))^2,$$

showing that the limiting variance is not sensitive to the underlying density but depends only on the dimension.

Still in the setting of general $\kappa$, the inequality $\int_A \kappa(x)^{(d-2\gamma)/d}\, dx \geq (\int_A \kappa(x)^{(d-\gamma)/d}\, dx)^2$ implies that the right-hand side of (2.12) is strictly positive whenever $\Delta^{\xi}(1)$ is nondegenerate. Moreover, (2.12) implies that when $d = 2$, $\kappa$ a density on $A = [0,1]^2$, and $\gamma = 1$, which would be the case for total edge length functionals of graphs on vertex sets in $[0,1]^2$, then the limiting variance of $n^{-1/2}\bar{\rho}^{\xi}_{n}(A)$ equals $V^{\xi}(1) - (D^{\xi}(1))^2(\int_A \kappa(x)^{1/2}\, dx)^2$, which is *minimized when the underlying density $\kappa$ is uniform.*

(iii) A comparison of (2.4) and (2.8) shows that Poissonization contributes extra randomness which shows up in the limiting variance (2.8). To show nondegeneracy of $\Delta^{\xi}(\kappa(X))$, we need to appeal to the particular geometric structure of the underlying problem. This is done on a case by case basis and is already treated in many problems of interest [24]. The implicit finiteness of the right-hand side of (2.4) and (2.8) is made explicit in Section 4.3.

(iv) Our method of proof actually yields (2.4) whenever (2.2) holds for some $p > 2$ and $\xi$ is exponentially stabilizing. This modification requires a small modification of Lemma 4.2. Also, if $\xi$ satisfies stabilization (Definition 2.1) only when $0 < a < b$, then Theorems 2.1 and 2.2 hold provided that $\kappa$ is bounded away from zero.

(v) The condition $m \in [n/2, 3n/2]$ in (2.5) is needed in order to achieve an efficient coupling between Poissonized and de-Poissonized measures. See Lemma 6.2 for details.

(vi) Theorem 2.1 holds for arbitrary continuous $\kappa : \mathbb{R}^d \to \mathbb{R}^+$; that is, $\kappa$ need not be a probability density.

### 2.3. *Extensions of main results.*

2.3.1. *Random measures induced by marked point processes.* Theorems 2.1 and 2.2 extend to random measures induced by marked point processes. Let $(\mathcal{M}, \mathcal{F}, \nu)$ be a probability space of marks and let $\mathcal{P}_{\tau \times \nu}$ (resp. $\mathcal{P}_{\kappa \times \nu}$) be a Poisson point process on $\mathbb{R}^d \times \mathcal{M}$ with intensity measure $\tau \times \nu$ (resp. $\kappa \times \nu$). We say that $\xi$ stabilizes if Definition 2.1 holds with $B_R(x)$ replaced by $B_R^{\mathcal{M}}(x) := B_R(x) \times \mathcal{M}$, and $\mathcal{X}$ ranging over the finite subsets of $(\mathbb{R}^d \setminus B_R(x)) \times \mathcal{M}$. Write $\xi(x; \mathcal{P}_\tau)$ for $\xi(x; \mathcal{P}_{\tau \times \nu})$.

Let $X_i, i \geq 1$, be i.i.d. marked random variables with common law $d\kappa \times \nu$. $X_i'$ denotes the projection of $X_i$ on $\mathbb{R}^d$, $\mathcal{X}_n := \{X_1, \ldots, X_n\}$, and $\rho_n^\xi := \sum_{i=1}^n \xi_n(X_i; \mathcal{X}_n) \delta_{X_i'}$ the associated marked random measures on $\mathbb{R}^d$. Let $\Delta^\xi(\tau)$ denote the marked version of (2.6), that is,

$$\begin{align}
(2.13) \qquad \Delta^\xi(\tau) &:= \Delta^\xi(\tau \times \nu) \\
&= \Delta_0((\mathcal{P}_{\tau \times \nu} \cap B_S^{\mathcal{M}}(\mathbf{0})) \cup \mathcal{A})
\end{align}$$

for all finite $\mathcal{A} \subset (\mathbb{R}^d \setminus B_S(\mathbf{0})) \times \mathcal{M}$. Let $V^\xi(\tau) := V^\xi(\tau \times \nu)$ denote the marked version of $V^\xi(\tau)$.

The analog of Theorem 2.2 for marked processes is as follows:

THEOREM 2.3. *Let $\rho_n^\xi, n \geq 1$, denote the marked measures defined above. Then (2.7)–(2.10) hold with $D^\xi(\tau)$ and $V^\xi(\tau)$ replaced by their respective marked versions $\mathbb{E}[\Delta^\xi(\tau \times \nu)]$ and $V^\xi(\tau \times \nu)$.*

REMARKS. (i) Theorem 2.3 generalizes Theorem 3.1 of [25], which establishes a CLT for *functionals* of marked *homogeneous* samples.

(ii) Applications of Theorem 2.3 to random sequential packing and spatial birth–growth models are discussed in Section 3.2.

2.3.2. *Random measures induced by nontranslation-invariant functionals.* It takes just a little extra effort to use our general approach to prove CLTs for measures induced by nontranslation-invariant weights $\xi$. Although translation invariance is often present in the measures (i)–(iii), we envision situations where measures on $\mathbb{R}^d$ do not enjoy translation invariance, as would be the case if the metric on $\mathbb{R}^d$ changes from point to point.

Let $\xi(y; x, \mathcal{X}), y \in \mathbb{R}^d$, be a family of measurable $\mathbb{R}$-valued functions defined for all pairs $(x, \mathcal{X})$, where $\mathcal{X} \subset \mathbb{R}^d$ is finite and $x$ is an element of $\mathcal{X}$. In cases with $x \notin \mathcal{X}$, we abbreviate the notation $\xi(y; x, \mathcal{X} \cup \{x\})$ to $\xi(y; x, \mathcal{X})$. We assume for all $y$ that $\xi(y; x, \mathcal{X})$ is translation invariant in the pairs $(x, \mathcal{X})$, that is, for all $z \in \mathbb{R}^d$ and all pairs $(x, \mathcal{X})$,

$$\xi(y; x, \mathcal{X}) = \xi(y; x - z, \mathcal{X} - z).$$

$\xi(y; \cdot, \cdot)$ is a rule depending on $y \in \mathbb{R}^d$ which assigns a real value to all pairs $(x, \mathcal{X})$. We do *not* assume that $\xi$ is translation invariant in the triples $(y; x, \mathcal{X})$. We define $\xi(x; \mathcal{X}) := \xi(x; x, \mathcal{X})$ for all $x \in \mathbb{R}^d$ and for all $\lambda > 0$, we set

$$\xi_\lambda(y; x, \mathcal{X}) := \xi(y; \lambda^{1/d} x, \lambda^{1/d} \mathcal{X}) \quad \text{and} \quad \xi_\lambda(x; \mathcal{X}) := \xi(x; \lambda^{1/d} x, \lambda^{1/d} \mathcal{X}).$$

We will consider limit theorems for the random measures

$$\mu_{\lambda\kappa}^\xi := \sum_{x \in \mathcal{P}_{\lambda\kappa}} \xi_\lambda(x; x, \mathcal{P}_{\lambda\kappa}) \delta_x$$

and

$$(2.14) \qquad \rho_n^\xi := \sum_{i=1}^n \xi_n(X_i; X_i, \mathcal{X}_n) \delta_{X_i}.$$

$\xi$ is said to be *stabilizing* if $\xi(x; \mathcal{X}) = \xi(x; x, \mathcal{X})$ stabilizes in the sense of Definition 2.1. If the rules $\xi(y; \cdot, \cdot)$ are identical for all $y \in \mathbb{R}^d$, then $\xi(x - z; \mathcal{X} - z) = \xi(x; \mathcal{X})$ for all $z \in \mathbb{R}^d$, and we reduce to the translation-invariant setting of Theorems 2.1 and 2.2.

When $\xi$ is translation invariant, that is to say, when $\xi(x; x, \mathcal{X}) = \xi(x, \mathcal{X})$, stabilization guarantees that pair correlation functions for $\xi$ decay suitably fast enough with respect to the interpoint distance. It also guarantees that the pair correlation function with respect to nonhomogeneous samples can be closely approximated by the pair correlation function with respect to homogeneous samples. However, for nontranslation-invariant $\xi$, a suitable approximation of pair correlation functions is not possible without some continuity of $\xi$ with respect to its first argument. This motivates the following definition.

DEFINITION 2.4.   The function $\xi$ is *slowly varying in* $L^q$ (abbreviated $\xi \in$ SV$(q)$) if for all $\tau \in (0, \infty)$, any $x \in \mathbb{R}^d$, and any compact set $K$ containing $\mathbf{0}$:

$$\lim_{\lambda \to \infty} \sup_{y \in K} \mathbb{E}[|\xi_\lambda(x + \lambda^{-1/d} y; x, \mathcal{P}_{\lambda\tau}) - \xi_\lambda(x; x, \mathcal{P}_{\lambda\tau})|^q] = 0.$$

The following generalizes Theorem 2.1 to nontranslation-invariant $\xi$. For all $x \in \mathbb{R}^d$ and $\tau > 0$, let

$$V^\xi(x, \tau) := \mathbb{E}[\xi^2(x; x, \mathcal{P}_\tau)]$$
$$+ \int_{\mathbb{R}^d} (\mathbb{E}[\xi(x; \mathbf{0}, \mathcal{P}_\tau \cup y) \cdot \xi(x; y, \mathcal{P}_\tau \cup \mathbf{0})]$$
$$- \mathbb{E}[\xi(x; \mathbf{0}, \mathcal{P}_\tau)] \mathbb{E}[\xi(x; y, \mathcal{P}_\tau')]) \tau \, dy.$$

THEOREM 2.4.   (i) *If $\xi \in SV(\frac{4}{3})$ is stabilizing and satisfies (2.2) for some $p > 1$, then for all $f \in C(A)$,*

$$\lim_{\lambda \to \infty} \frac{\mathbb{E}[\langle f, \mu_{\lambda\kappa}^{\xi} \rangle]}{\lambda} = \int_A f(x) \mathbb{E}[\xi(x; \mathcal{P}_{\kappa(x)})] \kappa(x)\, dx, \tag{2.15}$$

*whereas if $\xi$ is polynomially stabilizing and satisfies (2.2) for $p = 4$, then for all $f \in C(A)$,*

$$\lim_{\lambda \to \infty} \frac{\mathrm{Var}[\langle f, \bar{\mu}_{\lambda\kappa}^{\xi} \rangle]}{\lambda} = \int_A f^2(x) V^{\xi}(x, \kappa(x)) \kappa(x)\, dx. \tag{2.16}$$

(ii) *If $\xi \in SV(\frac{4}{3})$ is exponentially stabilizing and satisfies (2.2) for all $p > 0$, then $\lambda^{-1/2} \bar{\mu}_{\lambda\kappa}^{\xi}$ converges as $\lambda \to \infty$ to a Gaussian field with covariance kernel $\int_A f_1(x) f_2(x) V^{\xi}(x, \kappa(x)) \kappa(x)\, dx$.*

Letting $H^{\xi}(\mathcal{X}) = \sum_{x \in \mathcal{X}} \xi(x; \mathcal{X})$, we say that $H$ is strongly stabilizing for all $\tau > 0$ and all $x \in \mathbb{R}^d$ if there exist a.s. finite random variables $S$ (a *radius of stabilization* of $H$) and $\Delta^{\xi}(\tau, x)$ such that with probability 1,

$$\Delta^{\xi}(\tau, x) = \Delta_x((\mathcal{P}_{\tau} \cap B_S(x)) \cup \mathcal{A}) \tag{2.17}$$

for all finite $\mathcal{A} \subset \mathbb{R}^d \setminus B_S(x)$.

It is easy to check that the measures (2.14) satisfy the law of large numbers in Theorem 2.2(i). The following de-Poissonized version of Theorem 2.4 generalizes Theorem 2.2 and shows that the normalized versions of the measures (2.14) converge to a Gaussian field as well. Put $D^{\xi}(x, \tau) := \mathbb{E}[\Delta^{\xi}(x, \tau)]$.

THEOREM 2.5.   *Let $\xi \in SV(\frac{4}{3})$. Assume that $H$ is strongly stabilizing and satisfies the bounded moments condition for $\kappa$. Then we have:*

(i) *If $\xi$ is polynomially stabilizing and satisfies (2.2) for $p = 4$, then for all $f \in C(A)$,*

$$\lim_{n \to \infty} \frac{\mathrm{Var}[\langle f, \bar{\rho}_n^{\xi} \rangle]}{n}$$
$$= \int_A f^2(x) V^{\xi}(x, \kappa(x)) \kappa(x)\, dx - \left( \int_A f(x) D^{\xi}(x, \kappa(x)) \kappa(x)\, dx \right)^2. \tag{2.18}$$

(ii) *If $\xi$ is exponentially stabilizing and satisfies (2.2) for all $p > 0$, then $n^{-1/2} \bar{\rho}_n^{\xi}$ converges as $n \to \infty$ to a Gaussian field with covariance kernel*

$$\int_A f_1(x) f_2(x) V^{\xi}(x, \kappa(x)) \kappa(x)\, dx$$
$$- \int_A f_1(x) D^{\xi}(x, \kappa(x)) \kappa(x)\, dx \int_A f_2(x) D^{\xi}(x, \kappa(x)) \kappa(x)\, dx. \tag{2.19}$$

(iii) *If the distribution of $\Delta^\xi(X, \kappa(X))$ is nondegenerate, then*

$$\lim_{n \to \infty} \frac{\mathrm{Var}[\bar\rho_n^\xi(A)]}{n} > 0.$$

REMARKS.   (i) Formulas (2.18) and (2.19) are, in general, difficult to evaluate explicitly. However, in the context of statistics involving one-dimensional spacings, these formulas are readily evaluated [6], thus extending existing CLTs for sum functions of spacings.

(ii) We have used the assumption $\xi \in SV(4/3)$ only for technical convenience and have not aimed to find the optimal choice of $SV(q)$. Higher moment assumptions on $\xi$ will in general require $\xi \in SV(q)$ for smaller values of $q$ (cf. Lemma 4.2).

2.3.3. *Random measures induced by graphs.*   Theorems 2.1 and 2.2 assume a special form when the random point measures are induced by graphs. We see this as follows. Let $\mathcal{X}$ be a locally finite point set and let $G := G(\mathcal{X})$ be a graph on $\mathcal{X}$. $G$ is translation invariant if translation by $y$ is a graph isomorphism from $G(\mathcal{X})$ to $G(y + \mathcal{X})$ for all $y \in \mathbb{R}^d$ and all locally finite $\mathcal{X}$. $G$ is scale invariant if scalar multiplication by $a$ induces a graph isomorphism from $G(\mathcal{X})$ to $G(a\mathcal{X})$ for all $\mathcal{X}$ and all $a > 0$. Given $G$ and a vertex $x \in \mathcal{X}$, let $\mathcal{E}(x; G(\mathcal{X}))$ be the set of edges incident to $x$ (or for the Voronoi graph, the set of edges whose planar duals in the Delaunay graph are incident to $x$), and let $|e|$ denote the length of an edge $e$.

For any $f \in \mathcal{F}(a, b)$, let $\mathcal{P}_{f,x}$ denote $\mathcal{P}_f$ together with a point at $x$.

DEFINITION 2.5.   $G$ stabilizes if for all $A \in \mathcal{A}'$, $0 \le a < b < \infty$, $\lambda > 0$, and $x \in \lambda A$, there exists an a.s. finite random variable $R(x) := R(x, \lambda, a, b, A)$ (a *radius of stabilization*) defined on $(\Omega, \mathcal{F}, P)$ such that for all $f \in \mathcal{F}(a, b)$, with $\mathrm{Supp}\, f = \lambda A$, and all finite $\mathcal{X} \subset \lambda A \setminus B_R(x)$, we have

$$\mathcal{E}(x; G(\mathcal{P}_{f,x} \cap B_R(x) \cup \mathcal{X})) = \mathcal{E}(x; G(\mathcal{P}_{f,x} \cap B_R(x))).$$

Given $\phi : \mathbb{R}^+ \to \mathbb{R}^+$, consider functionals of the type $\xi_\phi^G(x; \mathcal{X}) := \sum_{e \in \mathcal{E}(x; G(\mathcal{X}))} \phi(|e|)$; such functionals could represent, for example, the total length of $\phi$-weighted edges in $G$ incident to $x$, the number of edges in $G$ incident to $x$, or the number of edges in $G$ less than some specified length. These functionals induce the point measures

$$\mu_{\phi,\mathcal{X}}^G := \sum_{x \in \mathcal{X}} \sum_{e \in \mathcal{E}(x; G(\mathcal{X}))} \phi(|e|) \delta_x.$$

If $G$ is polynomially stabilizing (resp. exponentially stabilizing), then so is $\xi_\phi^G$ for any $\phi$. Given $p > 1$, say that $\xi_\phi^G$ is $L^p$ bounded if

$$(2.20) \qquad \sup_{\lambda > 0} \sup_{x \in \lambda^{1/d} A} \sup_{\mathcal{X} \in \mathcal{C}} \mathbb{E}\left[ \left( \sum_{e \in \mathcal{E}(x; G(\mathcal{P}_{\lambda\kappa} \cup \mathcal{X}))} \phi(\lambda^{1/d}|e|) \right)^p \right] < \infty.$$

Let $H_\phi^G(\mathcal{X})$ denote the total mass of $\mu_{\phi, \mathcal{X}}^G$. $H_\phi^G(\mathcal{X})$ is strongly stabilizing if for all $\tau > 0$ there exist a.s. finite random variables $S := S(\tau)$ and $\Delta_\phi^G(\tau)$ such that with probability 1

$$H_\phi^G(\mathcal{P}_\tau \cup \{\mathbf{0}\} \cap B_S(\mathbf{0}) \cup \mathcal{A}) - H_\phi^G(\mathcal{P}_\tau \cap B_S(\mathbf{0}) \cup \mathcal{A}) = \Delta_\phi^G(\tau)$$

for all finite $\mathcal{A} \subset \mathbb{R}^d \setminus B_R(\mathbf{0})$.

Let $\mathcal{X}_n := n^{1/d}(X_1, \ldots, X_n)$, with $X_i$ i.i.d. with density $\kappa$. Write

$$(2.21) \qquad \rho_{\phi, n}^G := \sum_{i=1}^n \sum_{e \in \mathcal{E}(n^{1/d} X_i; G(\mathcal{X}_n))} \phi(|e|) \delta_{X_i}.$$

The next result is the CLT counterpart to the main result of [26] and follows immediately from Theorem 2.2. There is obviously a Poisson version of Theorem 2.6, but we do not state it here. Put $D_\phi^G(\tau) := \mathbb{E}[\Delta_\phi^G(\tau)]$ and let $V_\phi^G$ be the function $V^\xi$ when $\xi := \xi_\phi^G$ is defined as above.

THEOREM 2.6. *Assume that the graph $G$ is translation and scale invariant. Let $X_i, i \geq 1$, be i.i.d. with density $\kappa$. Assume that $H_\phi^G$ is strongly stabilizing and satisfies the bounded moments condition.*

(i) *If $\xi_\phi^G$ satisfies* (2.20) *for $p = 4$, if $G$ is polynomially stabilizing, then for all $\tau > 0$,*

$$\lim_{\lambda \to \infty} \frac{\mathrm{Var}[H_\phi^G(\lambda^{1/d}(\mathcal{P}_{\lambda\tau} \cap [0,1]^d))]}{\lambda} = V_\phi^G(\tau) \cdot \tau$$

*and for all $f \in C(A)$,*

$$(2.22) \quad \begin{aligned} \lim_{n \to \infty} &\frac{\mathrm{Var}[\langle f, \rho_{\phi, n}^G \rangle]}{n} \\ &= \int_A f^2(x) V_\phi^G(\kappa(x)) \kappa(x) \, dx - \left( \int_A f(x) D_\phi^G(\kappa(x)) \kappa(x) \, dx \right)^2. \end{aligned}$$

(ii) *If $G$ is exponentially stabilizing, if $\xi_\phi^G$ satisfies* (2.20) *for all $p > 0$, then as $n \to \infty$, $n^{-1/2} \bar{\rho}_{\phi, n}^G$ converges to a Gaussian field with covariance kernel*

$$(2.23) \quad \begin{aligned} &\int_A f_1(x) f_2(x) V_\phi^G(\kappa(x)) \kappa(x) \, dx \\ &\quad - \int_A f_1(x) D_\phi^G(\kappa(x)) \kappa(x) \, dx \int_A f_2(x) D_\phi^G(\kappa(x)) \kappa(x) \, dx. \end{aligned}$$

REMARKS. (i) If $\phi(x) = x^p$, $p > 0$, then the integrals in (2.22) and (2.23) can be simplified using the identities $V_\phi^G(\tau) = V_\phi^G(1)\tau^{-2p/d}$ and $D_\phi^G(\tau) = D_\phi^G(1)\tau^{-p/d}$.

(ii) We may generalize Theorem 2.6 to treat nontranslation-invariant $\xi$. For example, let $\xi_\phi^G(x; x, \mathcal{X})$ be a functional which assigns to a point $x$ in the graph $G(\mathcal{X})$ a value which depends on the point $x \in \mathbb{R}^d$ (e.g., the value may depend upon the local metric structure at $x$). Such functionals are not translation invariant in the triples $(x; x, \mathcal{X})$. By applying an appropriate uniformization to curved surfaces, we can fit functionals on such surfaces into our set-up of nontranslation-invariant functionals of point processes on $\mathbb{R}^d$. This yields, for example, CLTs for functionals of graphs defined over curved surfaces, in particular functionals of Voronoi diagrams over surfaces [17].

**3. Applications.** Theorems 2.1–2.6 can be applied to point measures induced by random graphs, packing processes and germ–grain models. This extends previous results [2, 21, 24, 26], [5, 7, 8, 9, 11, 15, 25, 27] and [12, 13, 22, 26] to the weak limit setting as well as to the setting of interaction processes over nonhomogeneous point fields. We do not provide an encyclopedic treatment of applications and anticipate applications to other interaction processes on $\mathbb{R}^d$, including measures induced by continuum percolation models. The methods described here can be modified to extend and generalize the central limit theory for classical spacings and $\phi$-divergences; in this setting the functions $V^\xi$ and $\Delta^\xi$ may be determined explicitly, allowing us to compute the limiting variance explicitly as a function of the underlying density of points. We refer to [6] for complete details.

Throughout, we will often show the exponentially stabilizing condition by appealing to results of [24, 26], which involves a slightly stronger definition of stabilization.

3.1. *Random graphs.* We limit discussion to random graphs on $\mathbb{R}^d$ with the usual Euclidean metric, but since translation invariance of $\xi$ is not assumed, many results hold if the graphs are defined on curved spaces. Our discussion parallels that in [26]. We say that $\phi$ has *polynomial growth* if there exists $a < \infty$ such that $\phi(x) \leq C(1 + x^a)$ for all $x \in \mathbb{R}^+$.

3.1.1. *k-nearest neighbors graphs.* Let $k$ be a positive integer. Given a locally finite point set $\mathcal{X} \subset \mathbb{R}^d$, the $k$-nearest neighbors (undirected) graph on $\mathcal{X}$, denoted NG($\mathcal{X}$), is the graph with vertex set $\mathcal{X}$ obtained by including $\{x, y\}$ as an edge whenever $y$ is one of the $k$ nearest neighbors of $x$ and/or $x$ is one of the $k$ nearest neighbors of $y$. The $k$-nearest neighbors (directed) graph on $\mathcal{X}$, denoted NG$'(\mathcal{X})$, is the graph with vertex set $\mathcal{X}$ obtained by placing a directed edge between each point and its $k$ nearest neighbors. $k$-nearest neighbors graphs are translation and scale invariant. Given a binomial sample $X_1, \ldots, X_n$ of i.i.d. random variables with density $\kappa$, define the induced point measures $\rho_{\phi,n}^{\mathrm{NG}}$ and $\rho_{\phi,n}^{\mathrm{NG}'}$ as in (2.21).

THEOREM 3.1. *The random measures* $\rho_{\phi,n}^{\mathrm{NG}}$ *and* $\rho_{\phi,n}^{\mathrm{NG}'}$, $n \geq 1$, *satisfy* (2.22), (2.23) *if* $\phi$ *has polynomial growth,* $\mathrm{Supp}\,\kappa \in \mathcal{A}$, *and* $\kappa$ *is bounded away from infinity and zero on its support.*

If we set $\phi^G(|e|) = |e|/2$, then we obtain a CLT for the total edge length of the $k$-nearest neighbors graph on the nonhomogeneous point set $\mathcal{X}_n$ whenever $\mathrm{Supp}\,\kappa \in \mathcal{A}$ and $\kappa$ is bounded away from infinity and zero. This generalizes existing CLTs [2, 24] which only show CLTs for nearest neighbor graphs on homogeneous point sets. The convergence to a Gaussian limit (2.23) is new.

Still more generally, if $\phi^G(|e|) = |e|^p/2$, $p > 0$, then Theorem 3.1 yields a CLT for the $p$th *power-weighted* total edge length of the $k$-nearest neighbors graph on $\mathcal{X}_n$ when $\mathrm{Supp}\,\kappa \in \mathcal{A}$ and $\kappa$ is bounded away from infinity and zero. That is, there are constants $V^{\mathrm{NG}}(1)$ and $D^{\mathrm{NG}}(1)$ such that

$$
\begin{aligned}
(3.1) \quad & \lim_{n \to \infty} \frac{\mathrm{Var}[\langle f, \bar{\rho}_{\phi,n}^{\mathrm{NG}} \rangle]}{n} \\
& = V^{\mathrm{NG}}(1) \int_A f^2(x) \kappa(x)^{(d-2\gamma)/d}\, dx \\
& \quad - (D^{\mathrm{NG}}(1))^2 \left( \int_A f(x) \kappa(x)^{(d-\gamma)/d}\, dx \right)^2
\end{aligned}
$$

and $n^{-1/2} \bar{\rho}_{\phi,n}^{\mathrm{NG}}$ converges as $n \to \infty$ to a Gaussian field with covariance kernel

$$
\begin{aligned}
(3.2) \quad & V^{\mathrm{NG}}(1) \int_A f_1(x) f_2(x) \kappa(x)^{(d-2\gamma)/d}\, dx \\
& \quad - (D^{\mathrm{NG}}(1))^2 \int_A f_1(x) \kappa(x)^{(d-\gamma)/d}\, dx \int_A f_2(x) \kappa(x)^{(d-\gamma)/d}\, dx.
\end{aligned}
$$

Another application of Theorem 3.1 goes as follows. Fix $t > 0$. Let $\phi^G(|e|)$ be either 0 or 1 depending on whether the length $|e|$ of the edge $e$ is bounded by $t$ or not. Then (2.23) gives a CLT for the empirical distribution function of the rescaled lengths of the edges in the $k$-nearest neighbors graph on $\mathcal{X}_n$.

PROOF OF THEOREM 3.1. The proof is straightforward and essentially follows from existing arguments in [24] and [26]. For completeness we sketch the proof when $G(\mathcal{X})$ denotes $\mathrm{NG}(\mathcal{X})$; similar arguments apply when $G(\mathcal{X})$ denotes $\mathrm{NG}'(\mathcal{X})$. It will suffice to apply Theorem 2.6 and to show that NG stabilizes on elements of $\mathcal{F}(a, b)$ when $a > 0$ [recall Remark (iv) after Theorem 2.2]. Let $f \in \mathcal{F}(a, b)$ be arbitrary, where $0 < a < b < \infty$. As shown in Lemma 6.1 of [24] (even though the definition of stabilization there is slightly different), the set of edges incident to $x$ in $\mathrm{NG}(\mathcal{P}_{f,x})$ is unaffected by the addition or removal of points outside a ball of random almost surely finite radius $4R$, that is, the graph $G(\mathcal{X}) = \mathrm{NG}(\mathcal{X})$ is stabilizing. Moreover, $R$ is constructed as follows [24]. For each $t > 0$ construct six disjoint equilateral triangles $T_j(t), 1 \le j \le 6$, such that $x$ is a vertex of each triangle and such that each triangle has edge length $t$. Then $R$ is the minimum $t$ such that each triangle $T_j(t), 1 \le j \le 6$, contains at least $k + 1$ points from $\mathcal{P}_{f,x}$. Since $f$ is bounded away from zero, elementary properties of the Poisson point process give the desired exponential decay of $R$ and thus $4R$ decays exponentially as well. We verify the moments condition (2.20) as in the proof of Theorem 2.4 of [26]. Strong stabilization of $H$ is given by Lemma 6.1 of [24] and the bounded moments condition for $H$ is as in Lemma 6.2 of [24]. The positivity of the limiting variance is given by Lemma 6.3 of [24]. $\square$

### 3.1.2. Voronoi and Delaunay graphs.

Given a locally finite set $\mathcal{X} \subset \mathbb{R}^d$ and $x \in \mathcal{X}$, the locus of points closer to $x$ than to any other point in $\mathcal{X}$ is called the *Voronoi cell* centered at $x$. The graph on vertex set $\mathcal{X}$ in which each pair of adjacent cell centers is connected by an edge is called the *Delaunay graph* on $\mathcal{X}$; if $d = 2$, then the planar dual graph consisting of all boundaries of Voronoi cells is called the *Voronoi graph* generated by $\mathcal{X}$. Edges of the Voronoi graph can be finite or infinite. Let $\mathrm{DEL}(\mathcal{X})$ [resp. $\mathrm{VOR}(\mathcal{X})$] denote the collection of edges in the Delaunay graph (resp. the Voronoi graph) on $\mathcal{X}$. The Voronoi and Delaunay graphs are clearly scale and translation invariant. Define the induced point measures $\rho_{\phi,n}^{\mathrm{VOR}}$ and $\rho_{\phi,n}^{\mathrm{DEL}}$ as in (2.21).

THEOREM 3.2. *Let $d = 2$. The random measures $\rho_{\phi,n}^{\mathrm{VOR}}$ and $\rho_{\phi,n}^{\mathrm{DEL}}$ satisfy (2.22), (2.23) if $\phi$ has polynomial growth with $\phi(\infty) = 0$, $\mathrm{Supp}\,\kappa \in \mathcal{A}$ and if $\kappa$ is bounded away from infinity and zero on its support.*

The limits (2.22), (2.23) extend the results of Penrose and Yukich [24] and Avram and Bertsimas [2], who consider CLTs for the total edge length of Voronoi graphs over *homogeneous* samples. The convergence in distribution to a Gaussian limit (2.23) is new. Clearly, the analogs of (3.1) and (3.2) hold for the measures $\rho_{\phi,n}^{\mathrm{VOR}}$ and $\rho_{\phi,n}^{\mathrm{DEL}}$.

PROOF OF THEOREM 3.2. We will apply Theorem 2.6. The moments condition (2.20) is shown in Theorem 2.5 of [26]. We can verify as in [26] that $G(\mathcal{X})$ is stabilizing. Let $f \in \mathcal{F}(a,b)$ be arbitrary, with $0 < a < b < \infty$. We will show that the Voronoi cell centered at $x$ with respect to $\mathcal{P}_{f,x}$ is unaffected by changes beyond a random but a.s. finite distance $R$ from $x$. We only need to show $R$ has exponentially decreasing tails. This is done in a manner similar to that for the $k$-nearest neighbors graph. For each $t > 0$ construct 12 disjoint isosceles triangles $T_j(t), 1 \le j \le 12$, such that $x$ is a vertex of each triangle, such that each triangle has two edges of length $t$, where $T_j(t) \subset T_j(u)$ whenever $t < u$ and where $\bigcup_{t>0} \bigcup_{j=1}^{12} T_j(t) = \mathbb{R}^2$. If $R$ is the minimum $t$ such that each triangle $T_j(t), 1 \le j \le 12$, contains at least one point from $\mathcal{P}_{f,x}$, then $3R$ is a radius of stabilization (page 1037 of [24]). Since $f$ is bounded away from zero, elementary properties of the Poisson point process give the desired exponential decay of $3R$. We can verify the bounded moments condition on $H$ as in Lemma 8.1 of [24]. Strong stabilization of $H$ is proved in Section 8 of [24]. The positivity of the limiting variance is given by Lemma 8.2 of [24]. $\quad\square$

3.1.3. *Sphere of influence graph.* Given a locally finite set $\mathcal{X} \subset \mathbb{R}^d$, the sphere of influence graph $\mathrm{SIG}(\mathcal{X})$ is a graph with vertex set $\mathcal{X}$, constructed as follows: for each $x \in \mathcal{X}$ let $B(x)$ be a ball around $x$ with radius equal to $\min_{y \in \mathcal{X} \setminus \{x\}} \{|y - x|\}$. Then $B(x)$ is called the sphere of influence of $x$. Draw an edge between $x$ and $y$ iff the balls $B(x)$ and $B(y)$ overlap. The collection of such edges is the sphere of influence graph (SIG) on $\mathcal{X}$ and is denoted by $\mathrm{SIG}(\mathcal{X})$. It is clearly translation and scale invariant. Define the induced point measure $\rho_{\phi,n}^{\mathrm{SIG}}$ as in (2.21).

In Section 7 of [24], CLTs are proved for the total edge length, the number of components, and the number of vertices of fixed degree of SIG when the underlying sample is uniform. The following extends these results to nonuniform samples and also shows weak convergence of the associated measures. We also obtain a CLT and variance asymptotics for the total number of edges in the SIG on nonuniform samples, extending results of [14].

THEOREM 3.3. *The random measures $\rho_{\phi,n}^{\mathrm{SIG}}$ satisfy (2.22), (2.23) if $\phi$ has polynomial growth, $\mathrm{Supp}\,\kappa \in \mathcal{A}$ and $\kappa$ is bounded away from infinity and zero on its support.*

PROOF. We will apply Theorem 2.6 again. Let $f \in \mathcal{F}(a, b)$ be arbitrary, $0 < a < b < \infty$. As shown in [26], $G(\mathcal{P}_f)$ has moments of all orders and is stabilizing, so we only need to show exponential stabilization. However, this follows from the analysis of SIG in [24]. Consider an infinite cone $C$ with its vertex at $x$, subtending an angle of $\pi/6$ radians. As in [24], let $T$ be the distance from $x$ to its closest neighbor in $\mathcal{P}_{f,x} \cap C$, and if $Y$ is the point in $C \cap B_{6T}(x)$ closest to $x$, then note (page 1030 of [24]) that the configuration of points outside $B_{3|Y|}(x)$ has no effect on the set of points in $C$ connected to $x$. Thus, the radius of stabilization $R$ equals the maximum of $m$ i.i.d. copies of $3|Y|$, where $m$ is the minimum number of cones $C_1, \ldots, C_m$ congruent to $C$, each with vertex at $x$, whose union is $\mathbb{R}^d$. It is easy to check that $R$ has exponential tails. Thus $G(\mathcal{X})$ is exponentially stabilizing. The bounded moments condition on $H$ is as in Lemma 7.2 of [24] and strong stabilization of $H$ is as in Lemma 7.1 of [24]. The positivity of the limiting variance is given by Theorem 7.2 of [24]. □

### 3.2. Random packing.
We will use Theorem 2.3 to extend earlier results on random sequential packing [5, 7, 8, 9, 11, 15, 25, 27] to cases of nonhomogeneous input as well as to show the weak convergence of packing measures induced by Poisson and fixed input.

### 3.2.1. RSA packing.
The following prototypical random sequential packing model is of considerable scientific interest.

Let $B_{n,1}, B_{n,2}, \ldots, B_{n,n}$ be a sequence of $d$-dimensional balls of volume $n^{-1}$ whose centers are i.i.d. random $d$-vectors $X_1, \ldots, X_n$ with probability density function $\kappa : A \to [0, \infty)$. Without loss of generality, assume that the balls are sequenced in the order determined by marks (time coordinates) in $[0, 1]$. Let the first ball $B_{n,1}$ be *packed*, and recursively for $i = 2, 3, \ldots, N$, let the $i$th ball $B_{n,i}$ be packed iff $B_{n,i}$ does not overlap any ball in $B_{n,1}, \ldots, B_{n,i-1}$ which has already been packed. If not packed, the $i$th ball is discarded. The collection of centers of accepted balls induces a point measure on $A$, denoted $\mu_{n^{-1}}$. We call this the *random sequential packing measure* induced by balls (of volume $n^{-1}$) with centers arising from $\kappa$.

Packing models of this type arise in diverse disciplines, including physical, chemical and biological processes. In statistical mechanics, this model describes the irreversible deposition of colloidal particles or proteins onto a substrate. In this context, the model described above is known as the RSA model for hard spheres on a continuum substrate. When the ball centers belong to a stationary Poisson point process on $\mathbb{R}^d$, then this model is the Matérn hard-core process (page 163 of [29]). When the ball centers belong to a bounded region of $\mathbb{R}^d$, then this model is known in spatial statistics as the simple sequential inhibition model (page 308 of [30]).

The vast scientific literature on versions of RSA models (see [25] for references) contains an abundance of experimental results, but few rigorous mathematical results. In $d = 1$, Rényi [27] and Dvoretzky and Robbins [11] established LLNs and CLTs, respectively, for the total number of accepted balls. Coffman, Flatto, Jelenković and Poonen ([9], equation (2), Theorems 13 and 14) determine explicit formulae for some of the limiting constants in the LLN and CLT, but restrict attention to $d = 1$. In $d \geq 1$, Penrose and Yukich [25] establish the asymptotic normality of the number of accepted balls when the spatial distribution is uniform and also show [26] an LLN for the number of accepted balls when the spatial distribution is nonuniform. Baryshnikov and Yukich [5] establish weak convergence of the sequential packing measures in $d \geq 1$, but only for homogeneous Poisson input. Here we will use our general result for marked processes, Theorem 2.3, to establish convergence of the variance and also weak convergence in the case of nonhomogeneous input in $d \geq 1$. We will follow the set-up of [26].

For any finite point set $\mathcal{X} \subset A$, assume the points $x \in \mathcal{X}$ have time coordinates which are independent and uniformly distributed over the interval $[0, 1]$. Assume balls of volume $n^{-1}$ are centered at the points of $\mathcal{X}$ and arrive sequentially in an order determined by the time coordinates, and assume as before that each ball is packed or discarded according to whether or not it overlaps a previously packed ball. Let $\xi(x; \mathcal{X})$ be either 1 or 0 depending on whether the ball centered at $x$ is packed or discarded. Let $\xi_n(x; \mathcal{X}) = \xi(n^{1/d}x; n^{1/d}\mathcal{X})$, where $n^{1/d}x$ denotes scalar multiplication of $x$ and *not* the mark associated with $x$ and where balls centered at points of $n^{1/d}\mathcal{X}$ have volume 1. Let $H(\mathcal{X}) := \sum_{x \in \mathcal{X}} \xi(x; \mathcal{X})$ be the total number of balls packed. The random measure

$$\mu_n^\xi := \sum_{i=1}^n \xi_n(X_i; \{X_i\}_{i=1}^n)\delta_{X_i}$$

coincides with $\mu_{n^{-1}}$.

Straightforward modifications of [25, 4] show that $\xi$ is exponentially stabilizing. The strict positivity of $V^\xi(\tau)$ is shown in Theorem 1.2 of [25]. Since $\xi$ is bounded it satisfies the moments condition (2.2). By Section 5 of [25], $H$ satisfies the bounded moments condition and strong stabilization. Therefore, Theorem 2.3 yields the following CLT.

THEOREM 3.4. *The random measures $\mu_n^\xi, n \geq 1$, satisfy* (2.8) *and* (2.9).

Theorem 3.4 shows asymptotic normality of the total number of accepted balls and generalizes [4, 25] to the case of nonhomogeneous input.

3.2.2. *Spatial birth–growth models.* Consider the following spatial birth–growth model in $\mathbb{R}^d$. Seeds are born at random locations $X_i \in \mathbb{R}^d$ at times $T_i$, $i = 1, 2, \ldots$, according to a unit intensity homogeneous spatial temporal Poisson point process $\Psi := \{(X_i, T_i) \in \mathbb{R}^d \times [0, \infty)\}$. When a seed is born, it forms a cell by growing radially in all directions with a constant speed $v \geq 0$. Whenever one growing cell touches another, it stops growing in that direction. Initially the seed takes the form of a ball of radius $\rho_i \geq 0$ centered at $X_i$. If a seed appears at $X_i$ and if the ball centered at $X_i$ with radius $\rho_i$ overlaps any of the existing cells, then the seed is discarded.

In the special case when the growth rate $v = 0$ and $\rho_i$ is constant, this model reduces to the RSA packing model. In the alternative special case where all initial radii are zero a.s., the model is known as the Johnson–Mehl model, originally studied in model crystal growth, and is described in Stoyan, Kendall and Mecke [29]. Chiu and Quine [8] show that the number of seeds accepted inside a cube $Q_\lambda$ of volume $\lambda$ by time $t$ satisfies a CLT, but apart from numerical considerations, their arguments do not preclude the possibility of limiting normal random variable with zero variance [7]. Penrose and Yukich [25] consider a modification of this model in which all seeds outside $Q_\lambda$ are automatically rejected, while the rules for seeds inside $Q_\lambda \times [0, \infty)$ are as above. They establish a CLT for this model and show that the limiting variance is strictly positive (page 295 of [25]), thus implying that the CLT of [8] is nondegenerate.

If seeds are born at random locations $X_i \in A$, it is natural to study the *spatial distribution* of accepted seeds. As far as we know, this problem has not been investigated. We may use Theorem 2.3 to establish the weak convergence of the random measure induced by the locations of the accepted seeds.

For any finite point set $\mathcal{X} \subset A$, assume the points $x \in \mathcal{X}$ have i.i.d. marks over $[0, 1]$. A mark at $x \in \mathcal{X}$ represents the arrival time of a seed at $x$. Assume that the seeds are centered at the points of $\mathcal{X}$ and that they arrive sequentially in an order determined by the associated marks, and assume that each seed is accepted or rejected according to the rules above. Let $\xi(x; \mathcal{X})$ be either 1 or 0 according to whether the seed centered at $x$ is accepted or not. $H(\mathcal{X}) := \sum_{x \in \mathcal{X}} \xi(x; \mathcal{X})$ is the total number of seeds accepted and $\Delta^\xi(\tau)$ is as in (2.7).

As with RSA packing, let $X_1, \ldots, X_n$ be i.i.d. random variables with density $\kappa$ on $A$ and with marks in $[0, 1]$. The random measure

$$\sigma_n^\xi := \sum_{i=1}^n \xi_n(X_i; \{X_i\}_{i=1}^n) \delta_{X_i}$$

is the scaled spatial birth–growth measure on $A$ induced by $X_1, \ldots, X_n$. The next result, a consequence of Theorem 2.3, shows that the spatial birth–growth measures converge to a Gaussian field.

THEOREM 3.5. *The random measures $\sigma_n^\xi, n \geq 1$, satisfy* (2.8) *and* (2.9).

Theorem 3.5 generalizes [8] and extends [25] to the case of nonhomogeneous input.

3.2.3. *Related packing models.* (a) Theorem 3.4 extends to more general versions of the prototypical packing model. By following the stabilization analysis of [25], one can develop asymptotics in the finite input setting for the number of packed balls in the following general models: (i) models with balls replaced by particles of random size/shape/charge, (ii) cooperative sequential adsorption models and (iii) ballistic deposition models (see [25] for a complete description of these models). In each case, Theorem 2.3 yields weak convergence to a Gaussian limit of the random packing measures associated with the centers of the packed balls, whenever the balls have a density $\kappa : A \to [0, \infty)$.

(b) The above packing models describe convergence of measures arising as a result of dependently thinning a Poisson point process. Related ways of thinning processes include the *annihilating process*, described as follows. A clock is attached to each point (particle) in the process; when the clock for a chosen particle rings, then if the particle has itself not been annihilated, it annihilates its neighbors within a fixed radius. Clearly, once a particle is free from occupied neighboring sites, it remains there undisturbed and is fixed for all time. Thus in any finite region the process is unchanging after a finite time. This models the thinning of seedlings [31] and the resulting random point measure satisfies the CLT in Theorem 3.4.

### 3.3. *Germ–grain models.*

Germ–grain models form a central part of stochastic geometry and spatial statistics [12, 22]. Here we consider the limit theory of functionals and measures associated with germ–grain models. Such models fall within the scope of the general set-up of Heinrich and Molchanov [13], who were the first to develop a general limit theory for random measures induced by translation-invariant germ–grain models.

Let $T_i, i \geq 1$, be i.i.d. bounded random variables defined on $(\Omega, \mathcal{F}, P)$, independent of the i.i.d. random variables $X_i, i \geq 1$, which are also defined on $(\Omega, \mathcal{F}, P)$ and which have density $\kappa$. For simplicity, consider random grains having the representation $X_i + n^{-1/d} B_{T_i}(X_i)$ and consider the random set

$$\Xi_n := \bigcup_{i=1}^{n} (X_i + n^{-1/d} B_{T_i}(\mathbf{0})).$$

When the $X_i, i \geq 1$, are the realization of a Poisson point process, the set $\Xi_n$ is a scale-changed Boolean model in the sense of Hall [12], pages 141 and 233. Heinrich and Molchanov [13] exploit the translation invariance of such a model and establish a central limit theorem for the associated measures. For translation-invariant models, Heinrich and Molchanov [13] establish CLTs without assuming boundedness of $T_i$.

For all $u \in \mathbb{R}^d$, let $T(u)$ be a random variable with a distribution equal to that of $T_1$. For all $x \in \mathbb{R}^d$ and all point sets $\mathcal{X} \subset \mathbb{R}^d$, let $V(x, \mathcal{X})$ be the Voronoi cell around $x$ with respect to $\mathcal{X}$. Given $x \in \mathbb{R}^d$, let $L(x, \mathcal{X})$ denote the Lebesgue measure of the intersection of the random set $\bigcup_{u \in \mathcal{X}} B_{T(u)}(u)$ and $V(x, \mathcal{X})$.

The *volume measure* induced by $\Xi_n$ is

$$\mu_n^L := \sum_{i=1}^{n} L_n(X_i; \mathcal{X}_n) \delta_{X_i}$$

and the total volume of $n^{1/d}\Xi_n$ is given by $H^L(n^{1/d}\mathcal{X}_n) := \sum_{i=1}^{n} L_n(X_i; \mathcal{X}_n)$.

Since $T$ is bounded it follows that $L$ is exponentially stabilizing and that $H^L$ is strongly stabilizing. Moreover, since the functional $L$ is bounded by the volume of a Voronoi cell, it is clear from Section 3.1.2 that $L$ satisfies the moment condition (2.2) for all $p > 0$ and that $H^L$ satisfies the bounded moments condition (2.5) for any $\kappa \in \mathcal{F}_{a,b}, 0 < a \leq b < \infty$.

Therefore, for germ–grain models $\Xi_n$ given above we have thus proved:

THEOREM 3.6.   *Let the density $\kappa$ be bounded away from infinity and zero.*

(i) *For all $f \in C(A)$*

$$\begin{aligned}
(3.3) \quad & \lim_{n \to \infty} \frac{\mathrm{Var}[\langle f, \bar{\mu}_n^L \rangle]}{n} \\
& = \int_A f^2(x) V^L(x, \kappa(x)) \kappa(x) \, dx - \left( \int_A f(x) D^L(x, \kappa(x)) \kappa(x) \, dx \right)^2.
\end{aligned}$$

(ii) *As $n \to \infty$, $n^{-1/2}\bar{\mu}_n^L$ converges to a Gaussian field with covariance kernel*

$$\begin{aligned}
(3.4) \quad & \int_A f_1(x) f_2(x) V^L(x, \kappa(x)) \kappa(x) \, dx \\
& - \int_A f_1(x) D^L(x, \kappa(x)) \kappa(x) \, dx \int_A f_2(x) D^L(x, \kappa(x)) \kappa(x) \, dx.
\end{aligned}$$

REMARKS.   (i) We have confined attention to one of the simplest germ–grain models. Instead of balls $B_T$, one could assume that the grains have some distribution on the space of convex subsets of $\mathbb{R}^d$. We have also limited our discussion to volume functionals, but it should be clear that the approach above readily extends to other spatial statistics, including total curvature.

(ii) Theorem 3.6 shows that volume functionals satisfy a CLT over nonuniform point sets, adding to results of [13] and [12], Chapter 3.4, involving the vacancy functional for germ–grain models.

(iii) The LLN counterpart of Theorem 3.6 is established in [26] and is not reproduced here.

## 4. Proof of variance convergence (Poisson case).

4.1. *Correlation functions.*   The proof of (2.16) uses the objective method [1] together with correlation functions. To illustrate the method, it is instructive to first prove the limit (2.3). Recall that for all $f \in C(A)$, $\lambda^{-1}\mathbb{E}[\langle f, \mu_{\lambda\kappa}^{\xi} \rangle] = \int_A f(x)\mathbb{E}[\xi_\lambda(x; \mathcal{P}_{\lambda\kappa})]\kappa(x) \, dx$. The key observation lying at the heart of the objective method is that for any point $x \in A$ distant at least $(K/\lambda)^{1/d}$ from $\partial A$, $K$ large, $\xi_\lambda(x; \mathcal{P}_{\lambda\kappa})$ is well approximated by the candidate limiting random variable $\xi(x; \mathcal{P}_{\kappa(x)})$ in the sense that as $\lambda \to \infty$

$$
\begin{aligned}
(4.1) \qquad & |\mathbb{E}\xi_\lambda(x; \mathcal{P}_{\lambda\kappa}) - \mathbb{E}\xi(x; \mathcal{P}_{\kappa(x)})| \\
& \leq |\mathbb{E}\xi_\lambda(x; \mathcal{P}_{\lambda\kappa}) - \mathbb{E}\xi_\lambda(x; \mathcal{P}_{\lambda\kappa(x)})| \\
& \quad + |\mathbb{E}\xi_\lambda(x; \mathcal{P}_{\lambda\kappa(x)}) - \mathbb{E}\xi(x; \mathcal{P}_{\kappa(x)})| \to 0,
\end{aligned}
$$

where $\mathcal{P}_{\lambda\kappa(x)}$ is a Poisson point process on $\mathbb{R}^d$ with intensity $\lambda\kappa(x)$ coupled to $\mathcal{P}_{\lambda\kappa}$ as in the upcoming coupling (4.2).

Indeed, to prove (4.1), for any point $x \in A$ distant at least $(K/\lambda)^{1/d}$ from $\partial A$, consider the event $F_{K,\lambda}(x)$ that the radius of stabilization $R(\lambda^{1/d}x)$ at $\lambda^{1/d}x$ with respect to $\mathcal{P}_{\kappa(x)}$ is less than $K$ and that $\mathcal{P}_{\lambda\kappa} = \mathcal{P}_{\lambda\kappa(x)}$ on $B_{(K/\lambda)^{1/d}}(x)$. Then

$$
P[F_{K,\lambda}(x)^c] \leq P[R(\lambda^{1/d}x) > K] + \lambda \int_{B_{(K/\lambda)^{1/d}}(x)} |\kappa(y) - \kappa(x)| \, dy.
$$

By stabilization, we have $P[R(\lambda^{1/d}x) > K] \leq \varepsilon$ uniformly in $\lambda$ by choosing $K$ large enough. For such $K$, the Lebesgue point property of $x$ shows that the second term above can be made arbitrarily small for large $\lambda$ and thus $P[F_{K,\lambda}^c(x)] \leq 2\varepsilon$ for large $\lambda$. Bound $|\mathbb{E}\xi_\lambda(x; \mathcal{P}_{\lambda\kappa}) - \mathbb{E}\xi_\lambda(x; \mathcal{P}_{\lambda\kappa(x)})|$ by

$$
\begin{aligned}
& |\mathbb{E}[(\xi_\lambda(x; \mathcal{P}_{\lambda\kappa}) - \xi_\lambda(x; \mathcal{P}_{\lambda\kappa(x)})) \cdot \mathbf{1}_{F_{K,\lambda}(x)}]| \\
& \quad + |\mathbb{E}[(\xi_\lambda(x; \mathcal{P}_{\lambda\kappa}) - \xi_\lambda(x; \mathcal{P}_{\lambda\kappa(x)})) \cdot \mathbf{1}_{F_{K,\lambda}(x)^c}]|.
\end{aligned}
$$

The first term vanishes by the definition of $F_{K,\lambda}(x)$. The second term is bounded by a multiple of $\varepsilon$ by combining Hölder's inequality, the assumed $1 + \delta$ moment condition and $P[F_{K,\lambda}^c(x)] \leq 2\varepsilon$. Thus, for any point $x \in A$ distant at least $(K/\lambda)^{1/d}$ from $\partial A$, the first term on the right-hand side of (4.1) goes to zero and the second term goes to zero by stabilization.

The proof of the variance convergence (2.16) is more involved and requires some extra terminology. Let $\mathcal{P}'_{\lambda\kappa}$ be a Poisson point process equidistributed with and independent of $\mathcal{P}_{\lambda\kappa}$, that is, $\mathcal{P}'_{\lambda\kappa}$ is a copy of $\mathcal{P}_{\lambda\kappa}$. For all $\lambda \in \mathbb{R}^+$ and $x \in A$ we introduce two auxiliary homogeneous (independent) Poisson point processes $\widetilde{\mathcal{P}}_{\lambda\kappa(x)}$ and $\widetilde{\mathcal{P}}'_{\lambda\kappa(x)}$ defined on $(\Omega, \mathcal{F}, P)$ such that:

(i) $\widetilde{\mathcal{P}}_{\lambda\kappa(x)}$ and $\widetilde{\mathcal{P}}'_{\lambda\kappa(x)}$ have *constant* intensity on $A$ equal to $\lambda\kappa(x)$,

(ii) $\mathcal{P}_{\lambda\kappa}$ and $\widetilde{\mathcal{P}}_{\lambda\kappa(x)}$ are coupled in the sense that for any Borel subset $B \subset A$,

$$
(4.2) \qquad P[\mathcal{P}_{\lambda\kappa}(B) \neq \widetilde{\mathcal{P}}_{\lambda\kappa(x)}(B)] \leq \lambda \int_B |\kappa(y) - \kappa(x)| \, dy,
$$

and the same is true for $\mathcal{P}'_{\lambda\kappa}$ and $\widetilde{\mathcal{P}}'_{\lambda\kappa(x)}$.

The proof of the variance convergence (2.16) approximates the correlations of $\xi_\lambda(x; \mathcal{P}_{\lambda\kappa})$, $x \in \mathbb{R}^d$, by those of $\xi_\lambda(x; \mathcal{P}_{\lambda\kappa(x)})$, $x \in \mathbb{R}^d$. Thus, for all $x \in \mathbb{R}^d$ define

$$q_\lambda(x) := \mathbb{E}[\xi_\lambda^2(x; \mathcal{P}_{\lambda\kappa})] \quad \text{and} \quad \tilde{q}_\lambda(x) := \mathbb{E}[\xi_\lambda^2(x; \widetilde{\mathcal{P}}_{\lambda\kappa(x)})],$$

as well as the pair correlation function

$$c_\lambda(x, y) := \mathbb{E}[\xi_\lambda(x; x, \mathcal{P}_{\lambda\kappa} \cup y)\xi_\lambda(y; y, \mathcal{P}_{\lambda\kappa} \cup x)]$$
$$- \mathbb{E}[\xi_\lambda(x; x, \mathcal{P}_{\lambda\kappa})]\mathbb{E}[\xi_\lambda(y; y, \widetilde{\mathcal{P}}_{\lambda\kappa}')], \qquad x, y \in \mathbb{R}^d, x \neq y.$$

Abbreviating notation throughout and writing $\xi_\lambda(x; \mathcal{P}_{\lambda\kappa} \cup y)\xi_\lambda(y; \mathcal{P}_{\lambda\kappa} \cup x)$ for $\xi_\lambda(x; x, \mathcal{P}_{\lambda\kappa} \cup y)\xi_\lambda(y; y, \mathcal{P}_{\lambda\kappa} \cup x)$, we also have the pair correlation functions in the homogeneous intensity case:

$$\tilde{c}_\lambda(x, y) := \mathbb{E}[\xi_\lambda(x; \widetilde{\mathcal{P}}_{\lambda\kappa(x)} \cup y)\xi_\lambda(y; \widetilde{\mathcal{P}}_{\lambda\kappa(x)} \cup x)$$
$$- \xi_\lambda(x; \widetilde{\mathcal{P}}_{\lambda\kappa(x)})\xi_\lambda(y; \widetilde{\mathcal{P}}_{\lambda\kappa(x)}')], \qquad x \neq y,$$

and

$$\tilde{c}_\lambda^x(x, y) := \mathbb{E}[\xi_\lambda^x(x; \widetilde{\mathcal{P}}_{\lambda\kappa(x)} \cup y)\xi_\lambda^x(y; \widetilde{\mathcal{P}}_{\lambda\kappa(x)} \cup x)$$
$$- \xi_\lambda^x(x; \widetilde{\mathcal{P}}_{\lambda\kappa(x)})\xi_\lambda^x(y; \widetilde{\mathcal{P}}_{\lambda\kappa(x)}')], \qquad x \neq y.$$

Here we employ the notation $\xi_\lambda^x(z; \mathcal{X})$ for $\xi_\lambda(x; z, \mathcal{X})$. Clearly, the correlations $\tilde{c}_\lambda(x, y)$ and $\tilde{c}_\lambda^x(x, y)$ are *not* symmetric in $x$ and $y$, unlike $c_\lambda$. When $\lambda = 1$ we write simply $q(x)$ and $c(x, y)$ for $q_1(x)$ and $c_1(x, y)$, respectively, and similarly for $\tilde{q}, \tilde{c}$ and $\tilde{c}^x$. Denote the integral of $f \in C(A)$ with respect to a Borel measure $\mu$ on $\mathbb{R}^d$ by $\langle f, \mu \rangle$. Now

$$(4.3) \quad \lambda^{-1} \operatorname{Var}[\langle f, \mu_{\lambda\kappa}^\xi \rangle] = \lambda \langle f \otimes f, \mathbb{E}[\mu_{\lambda\kappa}^\xi \otimes \mu_{\lambda\kappa}^\xi - \mu_{\lambda\kappa}^\xi \otimes \mu_{\lambda\kappa}'^\xi] \rangle,$$

where $f \otimes f$ denotes the function $f(x)f(y)$ on the product $A \times A$, $\mu \otimes \nu$ stands for the product measure on $A \times A$ and $\mu_{\lambda\kappa}'^\xi$ is just an independent copy of $\mu_{\lambda\kappa}^\xi$. Considering the diagonal and off-diagonal terms, we may rewrite the integral (4.3) in terms of correlation functions

$$\lambda^{-1} \operatorname{Var}[\langle f, \mu_{\lambda\kappa}^\xi \rangle]$$
$$= \lambda \iint_{A \times A} f(y)f(x)c_\lambda(x, y)\kappa(x)\kappa(y) \, dx \, dy + \int_A f(x)^2 q_\lambda(x)\kappa(x) \, dx$$
$$= \int_A \kappa(x)f(x)\left[ f(x)q_\lambda(x) + \lambda \int_A f(y)c_\lambda(x, y)\kappa(y) \, dy \right] dx.$$

To show the desired asymptotics (2.16), we will first show for all $x \in A$ not too close to the boundary of $A$, that as $\lambda \to \infty$

$$(4.4) \quad \begin{aligned} &f(x)q_\lambda(x) + \lambda \int_A f(y)c_\lambda(x, y)\kappa(y) \, dy \\ &\quad - \left[ f(x)\tilde{q}(x) + \int_{\mathbb{R}^d} f(x)\tilde{c}^x(x, y)\kappa(x) \, dy \right] \to 0. \end{aligned}$$

Note that the bracketed expression in (4.4) is $f(x)V^\xi(x, \kappa(x))$.

4.2. *Properties of correlation functions.* Showing the limit (4.4) requires some properties of correlation functions. Using the definitions and the translation invariance of $\xi(y; x, \mathcal{X})$ in the pairs $(x, \mathcal{X})$, it is easy to verify that for all $x, y \in A$:

$$\tilde{q}_\lambda(x) = \tilde{q}(x) \quad \text{and} \quad \tilde{c}_\lambda^x(x, x + y) = \tilde{c}^x(x, x + \lambda^{1/d}y).$$

Also, if (2.2) holds for $p = 2$, then we have the following uniform bounds:

$$(4.5) \qquad \sup_{x,y \in \mathbb{R}^d, \lambda > 0} [\max[q_\lambda(x), \tilde{q}_\lambda(x), c_\lambda(x, y), \tilde{c}_\lambda(x, y), \tilde{c}_\lambda^x(x, y)]] < \infty.$$

Our next fact provides some crucial decay properties. Here and elsewhere $C$ denotes a constant whose value may change from line to line.

LEMMA 4.1. *Under the moment condition* (2.2) *with* $p = 4$, *we have*

$$[\max(|c_\lambda(x, y)|, |\tilde{c}_\lambda(x, y)|, |\tilde{c}_\lambda^x(x, y)|)] \leq C(r(\lambda^{1/d}|x - y|/2))^{1/2}.$$

PROOF. We prove only

$$|c_\lambda(x, y)| \leq C(r(\lambda^{1/d}|x - y|/2))^{1/2}$$

since the proof of the other two inequalities is identical. Let $R_x := R(\lambda^{1/d}x)$ and $R_y := R(\lambda^{1/d}y)$ be the radii of stabilization of $\xi$ for $\lambda^{1/d}x$ and $\lambda^{1/d}y$, respectively. Let $E := E_{x,y,\lambda}$ denote the event that $R(\lambda^{1/d}x)$ and $R(\lambda^{1/d}y)$ are both less than $\delta := \lambda^{1/d}|x - y|/2$ and note that $P[E^c] \leq Cr(\lambda^{1/d}|x - y|/2)$. On $E$ the stabilization balls $B_{R_x}(\lambda^{1/d}x)$ and $B_{R_y}(\lambda^{1/d}y)$ do not intersect and thus on $E$ we have $\xi_\lambda(x; \mathcal{P}_{\lambda\kappa} \cup y) = \xi_\lambda(x; \mathcal{P}_{\lambda\kappa})$ and $\xi_\lambda(y; \mathcal{P}_{\lambda\kappa} \cup x) = \xi_\lambda(y; \mathcal{P}_{\lambda\kappa})$ showing that

$$|\mathbb{E}[\xi_\lambda(x; \mathcal{P}_{\lambda\kappa} \cup y)\xi_\lambda(y; \mathcal{P}_{\lambda\kappa} \cup x)]$$
$$- \mathbb{E}[\xi_\lambda(x; \mathcal{P}_{\lambda\kappa} \cap B_\delta(\lambda^{1/d}x))\xi_\lambda(y; \mathcal{P}_{\lambda\kappa} \cap B_\delta(\lambda^{1/d}y))]| \leq CP[E^c]$$

by Hölder's inequality and the moment condition (2.2) with $p = 4$. Using independence in the second expectation and the bound

$$|\mathbb{E}[\xi_\lambda(x; \mathcal{P}_{\lambda\kappa} \cap B_\delta(\lambda^{1/d}x))] - \mathbb{E}[\xi_\lambda(x; \mathcal{P}_{\lambda\kappa})]| \leq CP[E^c].$$

we are done. □

The next lemma shows that $q_\lambda$ and $c_\lambda$ are closely approximated by their "uniform versions" $\tilde{q}_\lambda$ and $\tilde{c}_\lambda^x$, respectively. Compactness of $A$ and the continuity of $\kappa$ and $f$ imply uniform continuity, so we fix *moduli of continuity* $t_\kappa, t_f : \mathbb{R}^+ \to \mathbb{R}^+$ such that for any $x, y \in A: |x - y| \leq \delta$, $|\kappa(x) - \kappa(y)| \leq t_\kappa(\delta)$ and $|f(x) - f(y)| \leq t_f(\delta)$.

LEMMA 4.2. *Assume that $\xi \in SV(\frac{4}{3})$ satisfies the moment condition (2.2) for $p = 4$ and is polynomially stabilizing. Then there exists a function $e \colon \mathbb{R}^+ \to \mathbb{R}^+$, decreasing to 0, and a function $\delta \colon \mathbb{R}^+ \to \mathbb{R}^+$, increasing to $\infty$, such that $\delta/\lambda \to 0$ and*

(i) *$\forall x \in A$ distant at least $(\delta/\lambda)^{1/d}$ from $\partial A$, $|q_\lambda(x) - \tilde{q}_\lambda(x)| \le e(\lambda)$,*

(ii) *$\forall x, y \in A$, each distant at least $(\delta/\lambda)^{1/d}$ from $\partial A$,*

$$\delta(\lambda)|c_\lambda(x,y) - \tilde{c}_\lambda^x(x,y)| \le e(\lambda),$$

(iii) *as $\lambda \to \infty$, the function $\delta$ satisfies*

$$\delta(\lambda)t_f((\delta(\lambda)/\lambda)^{1/d}) \to 0; \qquad \delta(\lambda)t_\kappa((\delta(\lambda)/\lambda)^{1/d}) \to 0.$$

PROOF. It is clear that one can always find a function $\delta := \delta(\lambda) \to \infty$ as $\lambda \to \infty$ such that $\delta/\lambda \to 0$ and implication (iii) holds, and even more strongly, that

$$\text{(4.6)} \qquad \delta(\lambda)^2 t_\kappa((\delta(\lambda)/\lambda)^{1/d}) \to 0 \qquad \text{as } \lambda \to \infty.$$

Fix such $\delta$.

For any $\lambda > 0$ and $x \in A$, $x$ distant at least $(\delta/\lambda)^{1/d}$ from $\partial A$, consider the event $\Omega_x := \Omega_{x,\lambda,\delta}$ that the radius of stabilization $R(\lambda^{1/d}x)$ of $\xi$ is less than $\delta$, and that $\mathcal{P}_{\lambda\kappa} = \widetilde{\mathcal{P}}_{\lambda\kappa(x)}$ on $B_{(\delta/\lambda)^{1/d}}(x)$, that is, $\lambda^{1/d}\mathcal{P}_{\lambda\kappa} = \lambda^{1/d}\widetilde{\mathcal{P}}_{\lambda\kappa(x)}$ on $B_{\delta^{1/d}}(\lambda^{1/d}x)$. By polynomial stabilization, by definition of $t_\kappa$, as well as by the coupling estimate (4.2), the probability of the complement of $\Omega_x$ is

$$\text{(4.7)} \qquad P[\Omega_x^c] \le \omega_d \delta t_\kappa((\delta/\lambda)^{1/d}) + a_1 \delta^{-a_2},$$

where $a_1 > 0$ and $a_2 > 2$ are constants and $\omega_d$ is the volume of the unit ball in $\mathbb{R}^d$. Thus $P[\Omega_x^c]$ tends to 0 as $\lambda \to \infty$ by (4.6). To show (i) we need to show $|\mathbb{E}[\xi_\lambda^2(x; \mathcal{P}_{\lambda\kappa}) - \xi_\lambda^2(x; \mathcal{P}_{\lambda\kappa(x)})]| \to 0$ as $\lambda \to \infty$:

$$|\mathbb{E}[\xi_\lambda^2(x; \mathcal{P}_{\lambda\kappa}) - \xi_\lambda^2(x; \mathcal{P}_{\lambda\kappa(x)})]|$$
$$\le \mathbb{E}[|\xi_\lambda^2(x; \mathcal{P}_{\lambda\kappa}) - \xi_\lambda^2(x; \mathcal{P}_{\lambda\kappa(x)})| \mathbf{1}_{\Omega_x}]$$
$$+ \mathbb{E}[|\xi_\lambda^2(x; \mathcal{P}_{\lambda\kappa}) - \xi_\lambda^2(x; \mathcal{P}_{\lambda\kappa(x)})| \mathbf{1}_{\Omega_x^c}].$$

The first term vanishes by the definition of $\Omega_x$ and the definition of $\xi_\lambda(x; \mathcal{P}_{\lambda\kappa(x)}) = \xi(\lambda^{1/d}x; \mathcal{P}_{\kappa(x)})$. Hölder's inequality, the moment condition (2.2) with $p = 4$, and the bound (4.7) show that the second term vanishes as $\lambda \to \infty$. This proves (i).

For the proof of (ii) it suffices to show that there exists a function $\delta(\lambda) \to \infty$ such that both

$$\text{(4.8)} \qquad \delta(\lambda)|c_\lambda(x,y) - \tilde{c}_\lambda(x,y)| \le e(\lambda)$$

and

$$(4.9) \qquad \delta(\lambda)|\tilde{c}_\lambda(x,y) - \tilde{c}_\lambda^x(x,y)| \le e(\lambda).$$

We first show the bound (4.8). By Lemma 4.1 it is enough to show for all $|x - y| < \delta/\lambda^{1/d}$ that $\delta(\lambda)|\mathbb{E}[\xi_\lambda(x; \mathcal{P}_{\lambda\kappa} \cup y)\xi_\lambda(y; \mathcal{P}_{\lambda\kappa} \cup x) - \xi_\lambda(x; \widetilde{\mathcal{P}}_{\lambda\kappa(x)} \cup y)\xi_\lambda(y; \widetilde{\mathcal{P}}_{\lambda\kappa(x)} \cup x)]| \to 0$ as $\lambda \to \infty$.

We proceed as in the proof of (i), but now consider the event $\Omega_{x,y}$ that the radii of stabilization $R(\lambda^{1/d}x)$ and $R(\lambda^{1/d}y)$ of $\xi$ for $\lambda^{1/d}x$ and $\lambda^{1/d}y$, respectively, are both less than $\delta$, that $\mathcal{P}_{\lambda\kappa} = \widetilde{\mathcal{P}}_{\lambda\kappa(x)}$ on the ball $B_{(\delta/\lambda)^{1/d}}(x)$, and that $\mathcal{P}'_{\lambda\kappa} = \widetilde{\mathcal{P}}'_{\lambda\kappa(x)}$ on the ball $B_{(\delta/\lambda)^{1/d}}(y)$. Since $y$ is within $\delta/\lambda^{1/d}$ of $x$, the probability that $\mathcal{P}_{\lambda\kappa} \ne \widetilde{\mathcal{P}}_{\lambda\kappa(x)}$ on the ball $B_{(\delta/\lambda)^{1/d}}(y)$ is less than $\omega_d \delta t_\kappa((\delta/\lambda)^{1/d}) + a_1 \delta^{-a_2}$. Therefore, under polynomial stabilization

$$P[\Omega_{x,y}^c] \le 2\omega_d \delta t_\kappa((\delta/\lambda)^{1/d}) + 2a_1\delta^{-a_2}.$$

The triangle inequality, the moment condition (2.2) with $p = 4$ and Hölder's inequality give

$$\delta|\mathbb{E}[\xi_\lambda(x; \mathcal{P}_{\lambda\kappa} \cup y)\xi_\lambda(y; \mathcal{P}_{\lambda\kappa} \cup x) - \xi_\lambda(x; \widetilde{\mathcal{P}}_{\lambda\kappa(x)} \cup y)\xi_\lambda(y; \widetilde{\mathcal{P}}_{\lambda\kappa(x)} \cup x)]|$$

$$\le a_3\delta(\mathbb{E}\xi^4)^{1/4}P[\Omega_{x,y}^c]^{3/4} \le a_4\delta(\delta t_\kappa((\delta/\lambda)^{1/d}) + \delta^{-a_2})^{3/4}.$$

We may similarly show

$$\delta|c_\lambda(x,y) - \tilde{c}_\lambda(x,y)| \le a_3\delta(\mathbb{E}\xi^4)^{1/4}P[\Omega_{x,y}^c]^{3/4} \le a_4\delta(\delta t_\kappa((\delta/\lambda)^{1/d}) + \delta^{-a_2})^{3/4},$$

which tends to zero as $\lambda \to \infty$ since $a_2 > 2$. Thus (4.8) is satisfied.

We now show the bound (4.9). Notice that the $SV(\frac{4}{3})$ assumption on $\xi$ implies automatically a somewhat stronger statement that the convergence is uniform not only on each fixed compact $K$, but also on the balls of radius $\delta_s(\lambda) \to \infty$ as $\lambda \to \infty$. Even more strongly, we have convergence to zero with rate $o((\delta_s(\lambda))^{-1})$ uniformly on balls of radius $\delta_s(\lambda)/\lambda^{1/d}$, that is,

$$\sup_{y\,:\,|x-y| \le \delta_s(\lambda)/\lambda^{1/d}} \delta_s(\lambda)\ \mathbb{E}[|\xi_\lambda(y; \widetilde{\mathcal{P}}_{\lambda\kappa(x)}) - \xi_\lambda^x(y; \widetilde{\mathcal{P}}_{\lambda\kappa(x)})|^{4/3}] \to 0$$

as $\lambda \to \infty$. Thus by Hölder's inequality ($p = 4, q = 4/3$) we have

$$\delta_s(\lambda)|\mathbb{E}[\xi_\lambda(x; \widetilde{\mathcal{P}}_{\lambda\kappa(x)})\xi_\lambda(y; \widetilde{\mathcal{P}}_{\lambda\kappa(x)}) - \xi_\lambda(x; \widetilde{\mathcal{P}}_{\lambda\kappa(x)})\xi_\lambda^x(y; \widetilde{\mathcal{P}}_{\lambda\kappa(x)})]|$$

$$\le \delta_s(\lambda)(\mathbb{E}[\xi_\lambda^4(x; \widetilde{\mathcal{P}}_{\lambda\kappa(x)})])^{1/4}(\mathbb{E}|\xi_\lambda(y; \widetilde{\mathcal{P}}_{\lambda\kappa(x)}) - \xi_\lambda^x(y; \widetilde{\mathcal{P}}_{\lambda\kappa(x)})|^{4/3})^{3/4} \to 0$$

uniformly over balls around $x$ of radius $\delta_s(\lambda)/\lambda^{1/d}$. Similarly,

$$\delta_s(\lambda)|\mathbb{E}[\xi_\lambda(x; \widetilde{\mathcal{P}}_{\lambda\kappa(x)})\xi_\lambda^x(y; \widetilde{\mathcal{P}}_{\lambda\kappa(x)}) - \xi_\lambda^x(x; \widetilde{\mathcal{P}}_{\lambda\kappa(x)})\xi_\lambda^x(y; \widetilde{\mathcal{P}}_{\lambda\kappa(x)})]| \to 0$$

as $\lambda \to \infty$. Therefore the bound (4.9) holds on such balls. On the complement of these balls, the bound (4.9) also holds by the polynomial decay of correlation functions given by Lemma 4.1. Therefore the bound (4.9) holds and the proof of Lemma 4.2 is complete.  $\square$

4.3. *Convergence of variance.* We establish the convergence of $\lambda^{-1} \operatorname{Var}[\langle f, \mu_{\lambda\kappa}^{\xi} \rangle]$ for all $f \in C(A)$. Convexity and compactness of $A$ implies the smoothness condition $\lim_{n\to\infty} n^{-1}\partial_r(n^{1/d}A) = 0$ for all $r > 0$, where $\partial_r(n^{1/d}A)$ denotes the volume of the $r$-neighborhood of the boundary of $n^{1/d}A$. Recalling (4.4), it suffices to show for all $x \in A$ distant at least $2(\delta/\lambda)^{1/d}$ from $\partial A$, that for large $\lambda$,

$$(4.10) \qquad f(x)q_\lambda(x) + \lambda \int_A f(y)c_\lambda(x,y)\kappa(y)\,dy$$

is close to

$$(4.11) \qquad f(x)\tilde{q}_\lambda(x) + \lambda \int_{\mathbb{R}^d} f(x)\tilde{c}_\lambda^x(x,y)\kappa(x)\,dy.$$

Without loss of generality, assume $\operatorname{Supp} f$ is the set $A$.

Lemma 4.2(i) implies that for all $x$ the difference of the first terms in (4.10) and (4.11) goes to zero as $\lambda \to \infty$. The difference of the integrals in (4.10) and (4.11) equals

$$(4.12) \qquad \lambda \int_{\mathbb{R}^d} [c_\lambda(x,y)f(y)\kappa(y) - \tilde{c}_\lambda^x(x,y)f(x)\kappa(x)]\,dy.$$

Let $\delta := \delta(\lambda)$ be as in Lemma 4.2 and let $B_{(\delta/\lambda)^{1/d}}(x)$ be the ball of radius $(\delta/\lambda)^{1/d}$ around $x$. To evaluate the integral (4.12), we integrate separately over $B_{(\delta/\lambda)^{1/d}}(x)$ and $\mathbb{R}^d \setminus B_{(\delta/\lambda)^{1/d}}(x)$. The integral over $B_{(\delta/\lambda)^{1/d}}(x)$ involves the difference

$$\lambda \int_{B_{(\delta/\lambda)^{1/d}}(x)} [c_\lambda(x,y)f(y)\kappa(y) - \tilde{c}_\lambda^x(x,y)f(x)\kappa(x)]\,dy,$$

which we split as

$$(4.13) \qquad \begin{aligned} &\lambda \int_{B_{(\delta/\lambda)^{1/d}}(x)} (c_\lambda(x,y) - \tilde{c}_\lambda^x(x,y))f(y)\kappa(y)\,dy \\ &+ \lambda \int_{B_{(\delta/\lambda)^{1/d}}(x)} \tilde{c}_\lambda^x(x,y)(f(y) - f(x))\kappa(y)\,dy \\ &+ \lambda \int_{B_{(\delta/\lambda)^{1/d}}(x)} \tilde{c}_\lambda^x(x,y)f(x)(\kappa(y) - \kappa(x))\,dy. \end{aligned}$$

The first integral is bounded by the product of $\lambda$, the volume of $B_{(\delta/\lambda)^{1/d}}(x)$ and the maximum of the integrand $(c_\lambda(x,y) - \tilde{c}_\lambda^x(x,y))f(y)\kappa(y)$. However, since $y$ is distant at least $(\delta/\lambda)^{1/d}$ from $\partial A$, the product goes to zero by Lemma 4.2(ii). The second and third integrals also tend to zero as $\lambda \to \infty$ by the bound (4.5) and Lemma 4.2(iii).

Since $f$ and $\kappa$ are bounded, the integral in (4.12) over $\mathbb{R}^d \setminus B_{(\delta/\lambda)^{1/d}}(x)$ is bounded by

$$(4.14) \qquad C \int_{\mathbb{R}^d \setminus B_{(\delta/\lambda)^{1/d}}(x)} [c_\lambda(x,y) + \tilde{c}_\lambda^x(x,y)]\,d(\lambda^{1/d}y),$$

which by Lemma 4.1 is bounded by

$$C \int_{\mathbb{R}^d \setminus B_{\delta^{1/d}}(x)} (r(|z-x|/2))^{1/2} + (r(|z-x|/2))^{1/2} \, dz.$$

The above integral is bounded by $2C\omega_d \int_{\delta^{1/d}/2}^{\infty} (r^\kappa(t))^{1/2} t^{d-1} \, dt$ which tends to zero as $\lambda \to \infty$ by assumption. We conclude that (4.14) converges to zero uniformly for all $x \in A$ distant at least $2(\delta/\lambda)^{1/d}$ from $\partial A$. Hence,

$$\lambda^{-1} \operatorname{Var}[\langle f, \mu_{\lambda\kappa}^\xi \rangle] - \int_A \kappa(x) f(x) \left[ f(x) \tilde{q}_\lambda(x) + \lambda \int_{\mathbb{R}^d} f(x) \tilde{c}_\lambda^x(x, x+y) \kappa(x) \, dy \right] dx$$

converges to zero as $\lambda \to \infty$. The equivalences $\tilde{q}_\lambda(x) = \tilde{q}(x)$ and $\tilde{c}_\lambda^x(x, x+y) = \tilde{c}^x(x, x + \lambda^{1/d} y)$ yield (4.4) as desired:

$$(4.15) \quad \begin{aligned} \lambda^{-1} &\operatorname{Var}[\langle f, \mu_\lambda^\xi \rangle] \\ &\to \int_A f(x) \left[ f(x) \tilde{q}(x) + \int_{\mathbb{R}^d} f(x) \tilde{c}^x(x, y) \kappa(x) \, dy \right] \kappa(x) \, dx. \end{aligned}$$

**5. Proof of Theorem 2.4.** We will only prove Theorems 2.4 and 2.5, since they are clearly a generalization of Theorems 2.1 and 2.2. We will first prove Theorem 2.4. We have already established Theorem 2.4(i) under the hypothesis that $\xi$ satisfies (2.2) for $p = 4$, and now to prove Theorem 2.4(ii) we assume that $\xi$ satisfies (2.2) for all $p > 0$.

5.1. *Cumulant measures.* Recall that $C(A)$ denotes the continuous functions $f : A \to \mathbb{R}$. To prove convergence of the finite-dimensional distributions of $\lambda^{-1/2} \bar{\mu}_{\lambda\kappa}^\xi, \lambda \geq 1$, in Theorem 2.4, it suffices to show for all test functions $f \in C(A)$ that the Laplace transform of the random variable $\lambda^{-1/2} \langle f, \bar{\mu}_{\lambda\kappa}^\xi \rangle$ converges as $\lambda \to \infty$ to the Laplace transform of a normal random variable with mean zero and variance $\frac{1}{2} \int_A f^2(x) V^\xi(x, \kappa(x)) \kappa(x) \, dx$. In other words, it suffices to show for all $f \in C(A)$ that

$$(5.1) \quad \lim_{\lambda \to \infty} \mathbb{E} \exp(\lambda^{-1/2} \langle -f, \bar{\mu}_{\lambda\kappa}^\xi \rangle) = \exp\left[ \frac{1}{2} \int_A f^2(x) V^\xi(x, \kappa(x)) \kappa(x) \, dx \right].$$

We will use the method of cumulants to show (5.1). We first recall the formal definition of cumulants. Let $W := \mathbb{R}^d$ and formally expand (5.1) in a power series in $f$ as follows:

$$(5.2) \quad \mathbb{E} \exp(\lambda^{-1/2} \langle -f, \bar{\mu}_{\lambda\kappa}^\xi \rangle) = 1 + \sum_{k=1}^{\infty} \frac{\lambda^{-k/2} \langle (-f)^k, M_\lambda^k \rangle}{k!},$$

where $f^k : \mathbb{R}^{dk} \to \mathbb{R}, k = 1, 2, \ldots$, is given by $f^k(v_1, \ldots, v_k) = f(v_1) \cdot \cdots \cdot f(v_k)$, and $v_i \in W, 1 \leq i \leq k$. $M_\lambda^k$ is a measure on $\mathbb{R}^{dk}$, *the $k$th moment measure* (page 130 of [10]).

We have

$$(5.3) \qquad dM_\lambda^k = m_\lambda(v_1, \ldots, v_k) \prod_{i=1}^k \kappa(v_i) \, d(\lambda v_i),$$

where the Radon–Nikodym derivative $m_\lambda(v_1, \ldots, v_k)$ is given by

$$(5.4) \qquad m_\lambda(v_1, \ldots, v_k) := \mathbb{E}\left[ \prod_{i=1}^k \bar{\xi}_\lambda(v_i; \mathcal{P}_{\lambda\kappa}) \right],$$

and where given $v_1, \ldots, v_k$ we abbreviate notation and write for all $1 \leq i \leq k$, $\bar{\xi}_\lambda(v_i; \mathcal{P}_{\lambda\kappa})$ for $\xi_\lambda(v_i; \mathcal{P}_{\lambda\kappa}) - \mathbb{E}[\xi_\lambda(v_i; \mathcal{P}_{\lambda\kappa})]$ and $\xi(v_i; \mathcal{P}_{\lambda\kappa})$ for $\xi(v_i; \mathcal{P}_{\lambda\kappa} \cup \{v_j\}_{j=1}^k)$. For each fixed $k$, the mixed moment on the right-hand side of (5.4) is finite uniformly in $\lambda$ by the moment bounds (2.2). Likewise, the $k$th summand in (5.2) is finite.

When the series (5.2) is convergent, the logarithm of the Laplace functional gives

$$\log\left[ 1 + \sum_{k=1}^\infty \frac{\lambda^{-k/2} \langle (-f)^k, M_\lambda^k \rangle}{k!} \right] = \sum_{l=1}^\infty \frac{\lambda^{-l/2} \langle (-f)^l, c_\lambda^l \rangle}{l!};$$

the signed measures $c_\lambda^l$ are *cumulant measures* (semi-invariants [20] or Ursell functions). Regardless of the validity of (5.2), all cumulants $c_\lambda^l, l = 1, 2, \ldots$, admit the representation

$$c_\lambda^l = \sum_{T_1, \ldots, T_p} (-1)^{p-1} (p-1)! M_\lambda^{T_1} \cdots M_\lambda^{T_p},$$

where $M_\lambda^{T_i}$ denotes a copy of the moment measure $M^{|T_i|}$ on the product space $W^{T_i}$ and where $T_1, \ldots, T_p$ ranges over all unordered partitions of the set $1, \ldots, l$ (see page 30 of [20]). More generally, $c_\lambda^T := c_\lambda(T)$ is the cumulant measure on $W^T$ with the representation

$$c_\lambda^T = \sum_{T_1, \ldots, T_p} (-1)^{p-1} (p-1)! M_\lambda^{T_1} \cdots M_\lambda^{T_p},$$

where $T_1, \ldots, T_p$ ranges over all unordered partitions of the set $T$. The first cumulant measure coincides with the expectation measure and the second cumulant measure coincides with the covariance measure. The cumulants $c_\lambda^l, l = 1, 2, \ldots$, all exist under the moment condition (2.2). In what follows we make critical use of the standard fact that if the cumulants $c_\lambda^l$ of a random variable $Z$ vanish for $l \geq 3$, then $Z$ has a normal distribution.

We will sometimes shorten notation and write $M^k, m$ and $c^l$ instead of $M_\lambda^k, m_\lambda$ and $c_\lambda^l$.

5.2. *Cluster measures.* Since $c_\lambda^1$ coincides with the expectation measure, we have $\langle f, c_\lambda^1 \rangle = 0$ for all $f \in C(A)$. We already know from Section 4 that $\lambda^{-1} \langle f^2, c_\lambda^2 \rangle = \lambda^{-1} \operatorname{Var}[\langle f, \mu_{\lambda\kappa}^\xi \rangle] \to \int_A f^2(x) V^\xi(x, \kappa(x)) \kappa(x) \, dx$. Thus, to prove (5.1), it will be enough to show for all $k \geq 3$ and all $f \in C(A)$ that $\lambda^{-k/2} \langle f^k, c_\lambda^k \rangle \to 0$ as $\lambda \to \infty$ (see, e.g., Lemma 3 of [28]).

A cluster measure $U_\lambda^{S,T}$ on $W^S \times W^T$ for nonempty $S, T \subset \{1, 2, \dots\}$ is defined by

$$U_\lambda^{S,T}(A \times B) = M_\lambda^{S \cup T}(A \times B) - M_\lambda^S(A) M_\lambda^T(B)$$

for all Borel $A$ and $B$ in $W^S$ and $W^T$, respectively.

Let $S_1$ and $S_2$ be a partition of $S$ and let $T_1$ and $T_2$ be a partition of $T$. A product of a cluster measure $U_\lambda^{S_1, T_1}$ on $W^{S_1} \times W^{T_1}$ with products of moment measures on $W^{S_2} \times W^{T_2}$ will be called a $(S, T)$ semi-cluster measure.

For each nontrivial partition $(S, T)$ of $\{1, \dots, k\}$ we next provide a representation of the $k$th cumulant $c^k$ as

$$c^k = \sum_{(S_1, T_1), (S_2, T_2)} \alpha((S_1, T_1), (S_2, T_2)) U^{S_1, T_1} M^{|S_2|} M^{|T_2|},$$

where the sum ranges over partitions of $(1, \dots, k)$ consisting of pairings $(S_1, T_1), (S_2, T_2)$, where $S_1, S_2 \subset S$ and $T_1, T_2 \subset T$, and where $\alpha((S_1, T_1), (S_2, T_2))$ are integer-valued prefactors.

In other words, for any nontrivial partition $(S, T)$ of $\{1, \dots, k\}$, we show that $c^k$ is a linear combination of $(S, T)$ semi-cluster measures. We were unable to find a proof of this in the literature and provide it.

LEMMA 5.1. *For each nontrivial partition $(S, T)$ of $\{1, \dots, k\}$ we have*

$$c^k = \sum_{(S_1, T_1), (S_2, T_2)} \alpha((S_1, T_1), (S_2, T_2)) U^{S_1, T_1} M^{|S_2|} M^{|T_2|}.$$

PROOF. The proof involves some notation and definitions. The moment measures $M^j$ are expressed in terms of the cumulants via

$$M^j = \sum_{T_1, \dots, T_p} c(T_1) \cdots c(T_p),$$

where the sum is over all partitions of $\{1, \dots, j\}$, that is, unordered collections $T_1, T_2, \dots, T_p$ of mutually disjoint subsets of $\{1, \dots, j\}$ whose union is $\{1, \dots, j\}$ (page 27 of [20], or [18]). Similarly, for any sets $S$ and $T$,

(5.5) $$M^{S \cup T} = \sum_{(S_1, T_1) \cdots (S_p, T_p)} c(S_1, T_1) \cdots c(S_p, T_p),$$

where the sum is over all partitions of $S \cup T$, where $S_i \subset S, T_i \subset T$. A typical element $(S_i, T_i), 1 \leq i \leq p$, of a partition thus involves a pair of sets, one a subset of $S$ and the other a subset of $T$. Some partitions of $S \cup T$ are such that the empty set appears in each pair $(S_i, T_i)$, $1 \leq i \leq p$. Call these the degenerate partitions.

We now prove Lemma 5.1. Split (5.5) as

$$(5.6) \qquad M^{S \cup T} = \sum_{\{\cdots\}} c(\cdots) \cdots c(\cdots) + \sum_{\{\cdots\}^*} c(\cdots) \cdots c(\cdots),$$

where $\{\cdots\}^*$ denotes degenerate partitions. The first sum contains the cumulant $c(S \cup T)$ as well as products of lower-order cumulants, that is, cumulants of the form $c(S_i \cup T_i)$, where $S_i \cup T_i$ is a proper subset of $\{1, \ldots, k\}$. Since each $c(S_i \cup T_i)$ is a product of moment measures, it follows that the first sum contains $c(S \cup T)$ as well as linear combinations of $(S, T)$ semi-cluster measures. The second sum is just the product of $M^S$ and $M^T$. Thus the cumulant measure $c(S \cup T)$ is

$$c(S \cup T) = M^{S \cup T} - M^S M^T + \text{l.c.,}$$

where l.c. denotes a linear combination of $(S, T)$ semi-cluster measures. Since $M^{S \cup T} - M^S M^T$ is a $(S, T)$ cluster measure, it follows that $c(S \cup T)$ is a linear combination of semi-clusters. In particular, if $(S, T)$ is a partition of $\{1, \ldots, k\}$, then $c^k$ is a linear combination of $(S, T)$ semi-cluster measures. $\square$

The following bound is critical for showing that $\lambda^{-k/2} \langle f, c_\lambda^k \rangle \to 0$ for $k \geq 3$ as $\lambda \to \infty$.

LEMMA 5.2. *If $\xi$ is exponentially stabilizing, then the functions $m_\lambda$ cluster exponentially; that is, for positive constants $A_{j,l}$ and $C_{j,l}$ one has uniformly*

$$|m_\lambda(x_1, \ldots, x_j, y_1, \ldots, y_l) - m_\lambda(x_1, \ldots, x_j) m_\lambda(y_1, \ldots, y_l)|$$

$$\leq A_{j,l} \exp(-C_{j,l} \delta \lambda^{1/d}),$$

*where $\delta := \min_{1 \leq i \leq j, 1 \leq p \leq l} |x_i - y_p|$ is the separation between the sets $(x_i)_{i=1}^j$ and $(y_p)_{p=1}^l$.*

PROOF. The proof is similar to that of Lemma 4.1. With probability at least $1 - \exp(-\delta \lambda^{1/d}/C)$, the radius of stabilization for each $\lambda^{1/d} x_i$, $1 \leq i \leq j$, and each $\lambda^{1/d} y_p$, $1 \leq p \leq l$, is less than $\lambda^{1/d} \delta$. Let $E_{j,l} := E_{j,l}(d)$ denote the event for which all such radii are less than $\lambda^{1/d} \delta$. On $E_{j,l}$ the stabilization balls do not intersect and therefore

$$\mathbb{E}\left[\prod_{i=1}^{j}\bar{\xi}_{\lambda}(x_i;\mathcal{P}_{\lambda\kappa})\prod_{p=1}^{l}\bar{\xi}_{\lambda}(y_p;\mathcal{P}_{\lambda\kappa})\mathbf{1}_{E_{j,l}}\right]$$

$$=\mathbb{E}\left[\prod_{i=1}^{j}\bar{\xi}_{\lambda}(x_i;\mathcal{P}_{\lambda\kappa})\mathbf{1}_{E_{j,l}}\right]\mathbb{E}\left[\prod_{p=1}^{l}\bar{\xi}_{\lambda}(y_p;\mathcal{P}_{\lambda\kappa})\mathbf{1}_{E_{j,l}}\right].$$

Hölder's inequality and the moment conditions imply there is a constant $A_{j,l}$ such that

$$|m_{\lambda}(x_1,\ldots,x_j,y_1,\ldots,y_l)-m_{\lambda}(x_1,\ldots,x_j)m_{\lambda}(y_1,\ldots,y_l)|\leq A_{j,l}(P[E_{j,l}^c])^{1/2}.$$

Since $P[E_{j,l}^c]$ decays exponentially in $\delta$, Lemma 5.2 follows.  □

The next lemma specifies decay rates for the cumulant measures. Such decay rates are useful in establishing moderate deviation principles and laws of the iterated logarithm for the measures $\bar{\mu}_{\lambda\kappa}^{\xi}$ [4]. Here we simply use the decay rates to conclude the proof of Theorem 2.4.

LEMMA 5.3. *For all $f\in C(A)$ and for all $k=2,3,\ldots$, we have $\lambda^{-k/2}\langle f^k, c_{\lambda}^k\rangle=O(\|f\|_{\infty}^k\lambda^{(2-k)/2})$.*

PROOF.   We need to estimate

$$\lambda^{-k/2}\int_{A^k}f(v_1)\cdots f(v_k)\,dc_{\lambda}^k(v_1,\ldots,v_k).$$

Let $\Delta_k$ denote the diagonal in $W^k$, that is, $v_1=v_2=\cdots=v_k$. For all $v:=(v_1,\ldots,v_k)\in A^k$, let $D_k(v)$ denote the distance to the diagonal.

Let $\Pi(k)$ be all partitions of $\{1,2,\ldots,k\}$ into two subsets $S$ and $T$. For all such partitions consider the subset $\sigma(S,T)$ of $A^S\times A^T$ having the property that $v\in\sigma(S,T)$ implies $d(x(v),y(v))\geq D_k(v)/k$, where $x(v)=v\cap A^S$ and $y(v)=v\cap A^T$. Since for every $v:=(v_1,\ldots,v_k)\in A^k$, there is a splitting $x:=x(v)$ and $y:=y(v)$ of $v$ such that $d(x,y)\geq D_k(v)/k$, it follows that $A^k$ is the union of the sets $\sigma(S,T),(S,T)\in\Pi(k)$. The key to the proof of Lemma 5.3 is to evaluate the cumulant $c_{\lambda}^k$ over each $\sigma(S,T)$. We then use Lemma 5.1 and adjust our choice of semi-clusters there to the particular choice of $(S,T)$.

By Lemma 5.1, the cumulant measure $dc_{\lambda}^k(v_1,\ldots,v_k)$ on $\sigma(S,T)$ is a linear combination of $(S,T)$ semi-cluster measures of the form

$$\sum_{(S_1,T_1),(S_2,T_2)}\alpha((S_1,T_1),(S_2,T_2))U^{S_1,T_1}M^{|S_2|}M^{|T_2|},$$

where the sum ranges over all partitions of $\{1, \ldots, k\}$ consisting of pairings $(S_1, T_1), (S_2, T_2)$, where $S_1, S_2 \subset S$ and $T_1, T_2 \subset T$, and where $\alpha((S_1, T_1), (S_2, T_2))$ are integer-valued prefactors.

Let $x$ and $y$ denote elements of $A^S$ and $A^T$, respectively. Let $\tilde{x}$ and $\tilde{y}$ denote elements of $A^{S_1}$ and $A^{T_1}$, respectively, and let $\tilde{x}^c$ denote the complement of $\tilde{x}$ with respect to $x$ and likewise with $\tilde{y}^c$. The integral of $f$ against an $(S, T)$ semi-cluster measure has the form

$$\lambda^{-k/2} \int_{\sigma(S,T)} f(v_1) \cdots f(v_k) \, d(M_\lambda^{|S_2|}(\tilde{x}^c) U_\lambda^{i+j}(\tilde{x}, \tilde{y}) M_\lambda^{|T_2|}(\tilde{y}^c)).$$

Letting $u_\lambda(\tilde{x}, \tilde{y}) := m_\lambda(\tilde{x}, \tilde{y}) - m_\lambda(\tilde{x}) m_\lambda(\tilde{y})$, and recalling (5.3), the above is bounded by

$$(5.7) \quad \lambda^{-k/2} \int_{\sigma(S,T)} f(v_1) \cdots f(v_k) m_\lambda(\tilde{x}^c) u_\lambda(\tilde{x}, \tilde{y}) m_\lambda(\tilde{y}^c) \prod_{i=1}^{k} \kappa(v_i) \, d(\lambda v_i).$$

Decompose the product measure $\prod_{i=1}^{k} \kappa(v_i) \, d(\lambda v_i)$ into two measures, one supported by the diagonal $\Delta_k$ and the other not. Off the diagonal, the integral (5.7) is bounded by

$$D\|f\|_\infty^k \lambda^{-k/2} \int_0^{\lambda^{1/d}} \exp(-Ct) P[D_k > t] \, dt = O(\lambda^{-k/2}\lambda),$$

since $u_\lambda$ decays exponentially with the distance to the diagonal (Lemma 5.2), the mixed moments $m_\lambda$ are uniformly bounded, and since the differential of a volume element of points at a distance greater than $t$ from the diagonal is bounded by the Lebesgue measure of the diagonal. Integrating over the diagonal measure $\lambda \kappa(v_1) \, dv_1$, and using the boundedness of $f$, we thus bound (5.7) by $D\|f\|_\infty^k \lambda^{-k/2}\lambda$ for some constant $D$. Since this estimate holds for all $\sigma(S, T), (S, T) \in \Pi(k)$, where $A^k$ is the finite union of sets $\sigma(S, T)$, Lemma 5.3 holds.   $\square$

## 6. De-Poissonization: proof of Theorem 2.5.
De-Poissonization involves a significant modification of de-Poissonization arguments for CLTs for translation-invariant functionals (Section 4 of [24]) defined over homogeneous point sets and thus we provide the details.

Before de-Poissonizing, we need the following definition. For any $\mathcal{X} \subset A$ and $f \in C(A)$, let

$$(6.1) \qquad\qquad H_n^f(\mathcal{X}) := \sum_{x \in \mathcal{X}} f(x) \xi_n(x; \mathcal{X}).$$

Letting $\mathcal{X}_m$ be a point process of $m$ i.i.d. random variables with density $\kappa$ on $A$, set $R_{m,n} := H_n^f(\mathcal{X}_{m+1}) - H_n^f(\mathcal{X}_m)$.

The following coupling lemma is inspired by and follows closely Lemma 4.2 in [24].

LEMMA 6.1. *Suppose $\xi$ is exponentially stabilizing for $\kappa$, and suppose $H$ is strongly stabilizing. Let $\varepsilon > 0$. Then there exists $\delta > 0$ and $n_0 \geq 1$ such that for all $n \geq n_0$ and all $m, m' \in [(1 - \delta)n, (1 + \delta)n]$ with $m < m'$, there exist random variables $X, X'$ with density $\kappa$ and a coupled family of variables $D := D(X), \ D' := D(X'), R := R(X, X'), R' := R'(X, X')$ with the following properties:*

(i) *$D$ and $D'$ each have the same distribution as $f(X)\Delta^{\xi}(X, \kappa(X))$;*
(ii) *$D$ and $D'$ are independent;*
(iii) *$(R, R')$ have the same joint distribution as $(R_{m,n}, R_{m',n})$;*
(iv) *$P[\{|D - R| > \varepsilon\} \cup \{|D' - R'| > \varepsilon\}] < \varepsilon$.*

PROOF. We will modify the proof of Lemma 4.2 of [24]. Suppose we are given $n$. Let $X, X', Y_1, Y_2, \ldots$ be i.i.d. random variables with density $\kappa$ on $A$. On the probability space $(\Omega, \mathcal{F}, P)$, let $\mathcal{P} = \mathcal{P}_{n\kappa}$ and $\mathcal{P}' = \mathcal{P}'_{n\kappa}$ be independent Poisson processes on $A$ with intensity measure $n\kappa(x)\,dx$.

Let $\mathcal{P}''$ be the point process consisting of those points of $\mathcal{P}$ which lie closer to $X$ than to $X'$ (in the Euclidean norm), together with those points of $\mathcal{P}'$ which lie closer to $X'$ than to $X$. Clearly $\mathcal{P}''$ is a Poisson process also having intensity measure $n\kappa(x)\,dx$ on $A$ and, moreover, it is independent of $X$ and of $X'$.

Let $N$ denote the number of points of $\mathcal{P}''$ (a Poisson variable with mean $n \cdot \text{vol}\, A$). Choose an ordering on the points of $\mathcal{P}''$, uniformly at random from all $N!$ possible such orderings. Use this ordering to list the points of $\mathcal{P}''$ as $W_1, W_2, \ldots, W_N$. Also, set $W_{N+1} = Y_1, W_{N+2} = Y_2, W_{N+3} = Y_3$ and so on.

Let
$$R := R(X, X') := H_n^f(\{W_1, \ldots, W_m, X\}) - H_n^f(\{W_1, \ldots, W_m\})$$
and
$$R' := R'(X, X') := H_n^f(\{W_1, \ldots, W_{m'-1}, X, X'\}) - H_n^f(\{W_1, \ldots, W_{m'-1}, X\}).$$
$X, X', W_1, W_2, W_3, \ldots$ are i.i.d. variables on $A$ with density $\kappa$, and therefore the pairs $(R, R')$ and $(R_{m,n}, R_{m',n})$ have the same joint distribution as claimed.

For all $x \in \mathbb{R}^d$ and $\tau > 0$, let $B(x, \tau)$ denote a stabilization ball for $\xi$ at $x$ with respect to $\mathcal{P}_\tau$. Recalling $\Delta_x(\mathcal{X}) := H(\mathcal{X} \cup x) - H(\mathcal{X})$ we put
$$D(x) := f(x)\Delta_x(\mathcal{P}_{\kappa(x)} \cap B(x, \kappa(x)))$$
and
$$D'(x') := f(x')\Delta_{x'}(\mathcal{P}'_{\kappa(x')} \cap B(x', \kappa(x'))).$$
Let $D := D(X), \ D' := D'(X')$. Then $D$ and $D'$ are independent, and by strong stabilization of $H$, given $(X, X') = (x, x')$, have the same distribution as $f(x)\Delta^\xi(x, \kappa(x))$ and $f(x')\Delta^\xi(x', \kappa(x'))$, respectively. It remains to show that $P[\{|D - R| > \varepsilon\} \cup \{|D' - R'| > \varepsilon\}] < \varepsilon$.

Without loss of generality, we may couple $\mathcal{P}_{\kappa(x)}$ and $n^{1/d}\mathcal{P}_{n\kappa}$ such that for all Borel sets $B \subset A$

$$P[\mathcal{P}_{n\kappa}(B) \neq \mathcal{P}_{n\kappa(x)}(B)] \leq n \int_B |\kappa(y) - \kappa(x)| \, dy.$$

Given $(X, X') = (x, x')$, for every $K > 0$ and $n = 1, 2, \ldots$, let $F_{K,n} := F_{K,n}(x)$ be the event that $B(n^{1/d}x, \kappa(x)) \subset B_K(n^{1/d}x)$ and that $\mathcal{P}_{\kappa(x)}$ and $n^{1/d}\mathcal{P}_{n\kappa}$ coincide on the ball $B_K(n^{1/d}x)$. As in the proof of the limit (2.3), we have by the uniform continuity of $\kappa$ and the Lebesgue point property of $x$ that for $K$ and $n$ large enough, $P[F_{K,n}^c] < \varepsilon/9$ uniformly in $x$. Similarly, let $F_{K,n} := F'_{K,n}(x')$ be the event that $B(n^{1/d}x', \kappa(x')) \subset B_K(n^{1/d}x')$ and that $\mathcal{P}_{\kappa(x')}$ and $n^{1/d}\mathcal{P}'_{n\kappa}$ coincide on $B_K(n^{1/d}x')$. For $K$ and $n$ large enough $P[(F'_{K,n})^c] < \varepsilon/9$ uniformly in $x'$.

Thus given $(X, X') = (x, x')$, on sets $F_{K,n}$ and $F'_{K,n}$ of probability at least $1 - \varepsilon/9$, we have

$$D(x) = f(x)\Delta_{n^{1/d}x}(n^{1/d}\mathcal{P}_{n\kappa} \cap B(n^{1/d}x, \kappa(x)))$$

and

$$D'(x') = f(x')\Delta_{n^{1/d}x'}(n^{1/d}\mathcal{P}'_{n\kappa} \cap B(n^{1/d}x', \kappa(x'))).$$

Thus given $(X, X') = (x, x')$, we need only show that

$$P[|f(x)\Delta_{n^{1/d}x}(n^{1/d}\mathcal{P}_{n\kappa} \cap B(n^{1/d}x, \kappa(x))) - R(x, x')| > \varepsilon] < \varepsilon$$

and

$$P[|f(x')\Delta_{n^{1/d}x'}(n^{1/d}\mathcal{P}_{n\kappa} \cap B(n^{1/d}x', \kappa(x'))) - R'(x, x')| > \varepsilon] < \varepsilon.$$

We will show the first bound only; the proof of the second bound is identical.

We now follow [24], page 1018. Choose $K$ large enough such that $P[S > K] < \varepsilon/9$. For all $w \in \mathbb{R}^d$ and $r > 0$, let $Q_r(w) = [-r, r]^d + w$ be the cube centered at $w$. Let $t_f$ denote the modulus of continuity of $f \in C(A)$ and find $b(n)$ such that $b(n)t_f(K/n^{1/d}) \to 0$. Given $\varepsilon$ and $K$ as above, let $n \geq n(\varepsilon, K)$ be so large that

(6.2)     $$\sup_x b(n) \cdot t_f(K/n^{1/d})\mathbb{E}|\xi(n^{1/d}W_1; n^{1/d}\mathcal{P}_{n\kappa} \cup n^{1/d}x)| < \varepsilon^2/72$$

and

(6.3)     $$b(n) \cdot t_f(K/n^{1/d})\mathbb{E}|\xi(n^{1/d}W_1; n^{1/d}\mathcal{P}_{n\kappa})| < \varepsilon^2/72.$$

Take $n$ so large that except on an event (denoted $E_0$) of probability less than $\varepsilon/9$, the positions of $n^{1/d}x$ and $n^{1/d}x'$ are sufficiently far from the boundary of $n^{1/d}A$ and from each other, that the cubes $Q_K(n^{1/d}x)$ and $Q_K(n^{1/d}x')$ are contained entirely within $n^{1/d}A$ (possible by the regularity of $\partial A$), and also are such that every point of $Q_K(n^{1/d}x)$ lies closer to $n^{1/d}x$ than to $n^{1/d}x'$ and every point of $Q_K(n^{1/d}x')$ lies closer to $n^{1/d}x'$ than to $n^{1/d}x$.

Set $\delta := \varepsilon(2K)^{-d}/18$. We assume $|m - n| \leq \delta n$ and $|m' - n| \leq \delta n$. For $n$ large enough, except on an event (denoted $E_1$) of probability at most $\varepsilon/9$, we have $|N - m| \leq 2\delta n = \varepsilon(2K)^{-d}n/9$, and likewise $|N - m'| \leq \varepsilon(2K)^{-d}n/9$.

Let $E$ be the event that the set of points of $n^{1/d}\{W_1, \ldots, W_m\}$ lying in $Q_K(n^{1/d}x)$ is not the same as the set of points of $\mathcal{P}$ lying in $Q_K(n^{1/d}x)$. This will happen either if one or more of the $(N - m)^+$ "discarded" points of $n^{1/d}\mathcal{P}''$ or one or more of the $(m - N)^+$ "added" points of $n^{1/d}\{Y_1, Y_2, \ldots\}$ lies in $Q_K(n^{1/d}x)$. For each added or discarded point, the probability of lying in $Q_K(n^{1/d}x)$ is at most $(2K)^d/n$, and so the probability of $E$, given that $E_1$ does not occur, is less than $\varepsilon/9$.

We now compute

$$
\begin{aligned}
P[&|D(x) - R(x,x')| > \varepsilon] \\
&\leq P[|D(x) - R(x,x')|\mathbf{1}_{F_{K,n} \cap E_0^c \cap E_1^c \cap E^c \cap \{S < K\}} > \varepsilon/2] \\
&\quad + P[|D(x) - R(x,x')|\mathbf{1}_{(F_{K,n} \cap E_0^c \cap E_1^c \cap E^c \cap \{S < K\})^c} > \varepsilon/2] \\
&\leq P[|D(x) - R(x,x')|\mathbf{1}_{F_{K,n} \cap E_0^c \cap E_1^c \cap E^c \cap \{S < K\}} > \varepsilon/2] \\
&\quad + P[F_{K,n}^c] + P[E_0] + P[E_1] + P[E \setminus E_1] + P[S > K].
\end{aligned}
$$

The last five terms are bounded by $\varepsilon/9$ for large $n$. Now consider the first probability. On the set $F_{K,n} \cap E_0^c \cap E_1^c \cap E^c \cap \{S < K\}$ the difference $|D(x) - R(x,x')|$ equals

$$
|f(x)\Delta_{n^{1/d}x}(n^{1/d}\mathcal{P}_{n\kappa} \cap B(n^{1/d}x, \kappa(x))) - R(x,x')|,
$$

which by strong stabilization of $H$ is bounded by (since $S < K$)

$$
\begin{aligned}
&\leq \sum_{n^{1/d}W_i \in Q_K(n^{1/d}x)} |f(W_i) - f(x)||\xi(n^{1/d}W_i; n^{1/d}\mathcal{P}_{n\kappa} \cup n^{1/d}x)| \\
&\quad + \sum_{n^{1/d}W_i \in Q_K(n^{1/d}x)} |f(W_i) - f(x)||\xi(n^{1/d}W_i; n^{1/d}\mathcal{P}_{n\kappa})|.
\end{aligned}
$$

By definition of $t_f$, the above is bounded by

$$
\begin{aligned}
&\leq t_f(K/n^{1/d}) \sum_{n^{1/d}W_i \in Q_K(n^{1/d}x)} |\xi(n^{1/d}W_i; n^{1/d}\mathcal{P}_{n\kappa} \cup n^{1/d}x)| \\
&\quad + t_f(K/n^{1/d}) \sum_{n^{1/d}W_i \in Q_K(n^{1/d}x)} |\xi(n^{1/d}W_i; n^{1/d}\mathcal{P}_{n\kappa})|.
\end{aligned}
$$

Let $N_K := \text{card}(n^{1/d}\{W_i\} \cap Q_K(n^{1/d}x))$. Then the first term in the above is bounded by

$$t_f(K/n^{1/d}) \sum_{n^{1/d}W_i \in Q_K(n^{1/d}x)} |\xi(n^{1/d}W_i; n^{1/d}\mathcal{P}_{n\kappa} \cup n^{1/d}x)|\mathbf{1}_{N_K \leq b(n)}$$

$$+ t_f(K/n^{1/d}) \sum_{n^{1/d}W_i \in Q_K(n^{1/d}x)} |\xi(n^{1/d}W_i; n^{1/d}\mathcal{P}_{n\kappa} \cup n^{1/d}x)|\mathbf{1}_{N_K > b(n)},$$

with a similar bound for the second term. Therefore, combining all of the above bounds

$$P[|D(x) - R(x,x')| > \varepsilon]$$
$$\leq 5\varepsilon/9 + P[|D(x) - R(x,x')|\mathbf{1}_{F_{K,n} \cap E_0^c \cap E_1^c \cap E^c \cap \{S < K\}} > \varepsilon/2]$$
$$(6.4) \quad \leq 5\varepsilon/9 + P\left[t_f(K/n^{1/d}) \sum_{i=1}^{b(n)} |\xi(n^{1/d}W_i; n^{1/d}\mathcal{P}_{n\kappa} \cup n^{1/d}x)| > \varepsilon/8\right]$$
$$+ P\left[t_f(K/n^{1/d}) \sum_{i=1}^{b(n)} |\xi(n^{1/d}W_i; n^{1/d}\mathcal{P}_{n\kappa})| > \varepsilon/8\right] + 2P[N > b(n)].$$

Using Chebyshev and the bounds (6.2) and (6.3), the second and third terms in (6.4) are bounded by $\varepsilon/9$ for $n$ large enough. For $n$ large, the last term in (6.4) is bounded by $\varepsilon/9$, since $N$ is a.s. finite.

Now integrate over all pairs $(x,x')$ to obtain the desired result. $\square$

The next lemma extends Lemma 4.3 of [24].

LEMMA 6.2. *Suppose $\xi$ is exponentially stabilizing and satisfies* (2.2) *for all $p > 0$. Suppose $H$ is strongly stabilizing and satisfies the bounded moments condition for $\kappa$. Let $(h(n))_{n \geq 1}$ be a sequence with $h(n)/n \to 0$ as $n \to \infty$. Then*

$$(6.5) \quad \lim_{n \to \infty} \sup_{n-h(n) \leq m \leq n+h(n)} |\mathbb{E}R_{m,n} - \mathbb{E}[f(X)\Delta^\xi(X, \kappa(X))]| = 0.$$

*Also*

$$(6.6) \quad \lim_{n \to \infty} \sup_{n-h(n) \leq m < m' \leq n+h(n)} |\mathbb{E}R_{m,n}R_{m',n} - \mathbb{E}[f(X)\Delta^\xi(X, \kappa(X))]^2| = 0$$

*and*

$$(6.7) \quad \lim_{n \to \infty} \sup_{n-h(n) \leq m \leq n+h(n)} |\mathbb{E}R_{m,n}^2| < \infty.$$

PROOF. We will follow the proof of Lemma 4.3 of [24]. Let $m$ be an arbitrary integer satisfying $n - h(n) \leq m \leq n + h(n)$. Let $\varepsilon > 0$. Provided $n$ is large enough, by Lemma 6.1 we can find coupled variables $D$ and $R$, with $D$ having the same distribution as $f(X)\Delta^\xi(X, \kappa(X))$, with $R$ having the same distribution as $R_{m,n}$, and with $P[|D - R| > \varepsilon] < \varepsilon$. Then

$$\mathbb{E} R_{m,n} = \mathbb{E} R = \mathbb{E}[D] + \mathbb{E}[(R - D)\mathbf{1}_{|R-D|>\varepsilon}] + \mathbb{E}[(R - D)\mathbf{1}_{|R-D|\leq\varepsilon}].$$

By Cauchy–Schwarz, the moments condition (2.2), and the fact that the bounded moment condition implies $E[D^2] < \infty$ (Lemma 4.1 of [24]), we have

$$\mathbb{E}|(R - D)\mathbf{1}_{|D-R|>\varepsilon}| \leq C\varepsilon^{1/2}.$$

Since $\varepsilon$ is arbitrarily small, (6.5) follows. The proof of (6.7) is similar and is omitted.

Next we consider $m, m'$ with $n - h(n) \leq m < m' \leq n + h(n)$. By Lemma 6.1, there are coupled variables $D, D', R, R'$ such that $D$ and $D'$ are independent and each has the same distribution as $f(X)\Delta^\xi(X, \kappa(X))$, $(R, R')$ have the same joint distribution as $(R_{m,n}, R_{m',n})$, and

$$P[\{|D - R| > \varepsilon\} \cup \{|D' - R'| > \varepsilon\}] < \varepsilon.$$

Now $\mathbb{E}[RR'] - \mathbb{E}[DD'] = \mathbb{E}[R(R' - D')] + \mathbb{E}[D'(R - D)]$. By Cauchy–Schwarz, we again obtain the bounds $\mathbb{E}[R|R' - D'|] < C\varepsilon^{1/2}$ and $\mathbb{E}[D'|R - D|] \leq C\varepsilon^{1/2}$. It follows that the difference $\mathbb{E}[RR'] - \mathbb{E}[DD']$ can be made arbitrarily small and (6.6) follows. $\quad\square$

PROOF OF THEOREM 2.5.   We first prove the limit (2.18). Given $f \in C(A)$, $\mathcal{X} \subset A$ and recalling (6.1), let $H_n^f := H_n^f(\mathcal{X}_n)$ and $H'_n{}^f := H_n^f(\mathcal{P}_{n\kappa})$. Assume that $\mathcal{P}_{n\kappa}$ is coupled to $\mathcal{X}_n$ by setting $\mathcal{P}_{n\kappa} := \{X_1, X_2, \ldots, X_{N_n}\}$, with $N_n$ an independent Poisson variable with mean $n$.

To prove (2.18), it is enough to show for all $f \in C(A)$ that

$$(6.8) \qquad \frac{H_n^f - \mathbb{E} H_n^f}{n^{1/2}} \to N(0, \tau_f^2),$$

where $N(0, \tau_f^2)$ denotes a mean zero normal random variable with variance

$$(6.9) \qquad \begin{aligned} \tau_f^2 := {} & \int_A f^2(x) V^\xi(x, \kappa(x)) \kappa(x)\, dx \\ & - \left( \int_A f(x) \mathbb{E}[\Delta^\xi(x, \kappa(x))] \kappa(x)\, dx \right)^2. \end{aligned}$$

Letting $\alpha := \mathbb{E}[f(X)\Delta^\xi(X, \kappa(X))]$, the first step is to prove that as $n \to \infty$,

$$(6.10) \qquad \mathbb{E}[n^{-1/2}(H'_n{}^f - H_n^f - (N_n - n)\alpha)^2] \to 0.$$

To do this, we employ the coupling Lemma 6.1 and follow pages 1019 and 1020 of [24] verbatim. The second step is to prove that

$$(6.11) \qquad \lim_{n\to\infty} \frac{\operatorname{Var} H_n^f}{n} = \tau_f^2.$$

However, this follows from the identity

$$n^{-1/2}H_n'^f = n^{-1/2}H_n^f + n^{-1/2}(N_n - n)\alpha + n^{-1/2}(H_n'^f - H_n^f - (N_n - n)\alpha).$$

The third term in the above has variance tending to zero by (6.10); the second term has variance $\alpha^2$ and is independent of the first term. Letting $\sigma_f^2 := \int_A f^2(x)V^\xi(x, \kappa(x))\kappa(x)\,dx$, it follows that

$$\sigma_f^2 = \lim_{n\to\infty}\frac{\mathrm{Var}(H_n'^f)}{n} = \lim_{n\to\infty}\frac{\mathrm{Var}(H_n^f)}{n} + \alpha^2,$$

that is, (6.11) holds. The limit (6.8) follows as on page 1020 of [24], thus establishing (2.18).

Let us now show (2.19). The above shows that the sequence of distributions $n^{-1/2}\langle\bar\rho_n^\xi, f\rangle$ tends to a limiting normal random variable $N(0, \tau_f^2)$ for every $f \in C(A)$. Taking $f = f_1 + f_2$ and using simple algebra shows that the limiting Gaussian field has the desired covariance matrix (2.19). This completes parts (i) and (ii) of Theorem 2.5.

To prove Theorem 2.5(iii), it suffices to show that

$$(6.12)\quad \begin{aligned}\lim_{n\to\infty}&\frac{\mathrm{Var}[H_n(\mathcal{X}_n)]}{n}\\&= \int_A V^\xi(x, \kappa(x))\kappa(x)\,dx - \left(\int_A \mathbb{E}[\Delta^\xi(x, \kappa(x))]\kappa(x)\,dx\right)^2 > 0.\end{aligned}$$

We accomplish this by modifying the approach in Section 5 of [24].

We write $H_n(\mathcal{X}_n) - \mathbb{E}H_n(\mathcal{X}_n)$ as a sum of martingale differences as follows. Let $\mathcal{F}_i = \sigma(X_{1,n}, \ldots, X_{i,n})$ and write $\mathbb{E}_i$ for conditional expectation given $\mathcal{F}_i$. Define martingale differences $D_i := \mathbb{E}_i H_n(\mathcal{X}_n) - \mathbb{E}_{i-1}H_n(\mathcal{X}_n)$. Then $H_n(\mathcal{X}_n) - \mathbb{E}H_n(\mathcal{X}_n) = \sum_{i=1}^n D_i$ and

$$\mathrm{Var}[H_n(\mathcal{X}_n)] = \sum_{i=1}^n \mathbb{E}[D_i^2].$$

It suffices to show that there exists a constant $C > 0$ such that for all $1 \le i \le n$, $\mathbb{E}[D_i^2] > C$.

Given $i \le m$, let $G_{i,m} = H_n(\mathcal{X}_m) - H_n(\mathcal{X}_m \setminus \{X_i\})$. Let $\widetilde{G}_{i,m} = H_n(\mathcal{X}_{m+1} \setminus \{X_i\}) - H_n(\mathcal{X}_m \setminus \{X_i\})$. Then $D_i = \mathbb{E}_i[G_{i,n} - \widetilde{G}_{i,n}]$. We set $\alpha := \mathbb{E}[\Delta^\xi(X, \kappa(X))]$ and using nondegeneracy, take $\delta > 0$ such that $P[\Delta^\xi(X, \kappa(X)) > \alpha + 4\delta] > 4\delta$.

Define $f : \mathbb{R} \to \mathbb{R}$ by $f(x) = 0$ for $x \le \alpha + \delta$ and $f(x) = 1$ for $x \ge \alpha + 2\delta$, interpolating linearly between $\alpha + \delta$ and $\alpha + 2\delta$. Let $Y_i := f(\mathbb{E}_i[G_{i,n}])$. The remainder of the proof consists in showing that for $n$ large and for $i \ge (1 - \varepsilon_3)n$, we have

$$\mathbb{E}[(G_{i,n} - \alpha)Y_i] \ge 4\delta^2 \quad\text{and}\quad \mathbb{E}[(\widetilde{G}_{i,n} - \alpha)Y_i] \le 2\delta^2.$$

These inequalities follow from Lemma 6.1 and pages 1021 and 1022 of [24]. Taken together, this implies for large $n$ and $i \geq (1 - \varepsilon_3)n$, that $\mathbb{E}[(G_{i,n} - \widetilde{G}_{i,n})Y_i] \geq 2\delta^2$. Using the fact that $Y_i$ is $\mathcal{F}_i$-measurable and lies in the range $[0, 1]$, we obtain

$$2\delta^2 \leq \mathbb{E}[Y_i \mathbb{E}_i(G_{i,n} - \widetilde{G}_{i,n})] \leq \mathbb{E}[|\mathbb{E}_i(G_{i,n} - \widetilde{G}_{i,n})|] = \mathbb{E}[|D_i|],$$

and hence, $\mathbb{E}[D_i^2] \geq [\mathbb{E}|D_i|]^2 \geq 4\delta^4 > 0$. Thus (6.12) holds, completing the proof of Theorem 2.5. $\square$

**Acknowledgments.** The authors deeply appreciate heplful comments of Mathew Penrose and Tomasz Schreiber who both made corrections on sections of earlier drafts of the paper. The authors gratefully thank an anonymous referee for copious comments leading to an improved exposition.

BELL LABORATORIES, LUCENT TECHNOLOGIES
600 MOUNTAIN AVENUE
MURRAY HILL, NEW JERSEY 07974-0636
USA
E-MAIL: ymb@research.bell-labs.com

DEPARTMENT OF MATHEMATICS
LEHIGH UNIVERSITY
BETHLEHEM, PENNSYLVANIA 18015
USA
E-MAIL: joseph.yukich@lehigh.edu